\documentclass[a4paper,twoside]{article}
\usepackage{a4}
\usepackage{amssymb}
\usepackage{amsmath}
\usepackage{upref}
\usepackage[pagebackref,colorlinks,citecolor=blue,linkcolor=blue,urlcolor=blue]{hyperref}
\usepackage[active]{srcltx}
\usepackage[dvipsnames]{color}
\allowdisplaybreaks[2] 
\newcount\minutes \newcount\hours
\hours=\time \divide\hours 60 \minutes=\hours
 \multiply\minutes-60
\advance\minutes \time
\newcommand{\klockan}{\the\hours:{\ifnum\minutes<10 0\fi}\the\minutes}
\newcommand{\tid}{\today\ \klockan}
\newcommand{\prtid}{\smash{\raise 10mm \hbox{\LaTeX ed \tid}}}
\renewcommand{\prtid}{}
\makeatletter  \headheight 10pt
\def\sectionmark#1{} 
\def\subsectionmark#1{}
\newcommand{\sectnr}{\ifnum \c@secnumdepth >\z@
                 \thesection.\hskip 1em\relax \fi}
\def\@evenhead{\footnotesize\rm\thepage\hfil\leftmark\hfil\llap{\prtid}}
\def\@oddhead{\footnotesize\rm\rlap{\prtid}\hfil\rightmark\hfil\thepage}
\def\tableofcontents{\section*{Contents} 
 \@starttoc{toc}}
\makeatother
\makeatletter
\def\@biblabel#1{#1.}
\makeatother
\makeatletter
\let\Thebibliography=\thebibliography
\renewcommand{\thebibliography}[1]{\def\@mkboth##1##2{}\Thebibliography{#1}
\addcontentsline{toc}{section}{References}
\frenchspacing 
\setlength{\@topsep}{0pt}
\setlength{\itemsep}{0pt}%
\setlength{\parskip}{0pt plus 2pt}%
}
\makeatother
\makeatletter
\def\mdots@{\mathinner.\nonscript\!.%
 \ifx\next,.\else\ifx\next;.\else\ifx\next..\else
 \nonscript\!\mathinner.\fi\fi\fi}
\let\ldots\mdots@
\let\cdots\mdots@
\let\dotso\mdots@
\let\dotsb\mdots@
\let\dotsm\mdots@
\let\dotsc\mdots@
\def\vdots{\vbox{\baselineskip2.8\p@ \lineskiplimit\z@
    \kern6\p@\hbox{.}\hbox{.}\hbox{.}\kern3\p@}}
\def\ddots{\mathinner{\mkern1mu\raise8.6\p@\vbox{\kern7\p@\hbox{.}}%
    \raise5.8\p@\hbox{.}\raise3\p@\hbox{.}\mkern1mu}}
\makeatother
\makeatletter
\let\Enumerate=\enumerate
\renewcommand{\enumerate}{\Enumerate%
\setlength{\@topsep}{0pt}
\setlength{\itemsep}{0pt}%
\setlength{\parskip}{0pt plus 1pt}%
\renewcommand{\theenumi}{\textup{(\alph{enumi})}}%
\renewcommand{\labelenumi}{\theenumi}%
}
\let\endEnumerate=\endenumerate
\renewcommand{\endenumerate}{\endEnumerate\unskip}
\makeatother
\makeatletter
\def\@seccntformat#1{\csname the#1\endcsname.\quad}
\makeatother
\newcommand{\authortitle}[2]{\author{#1}\title{#2}\markboth{#1}{#2}}
\newcommand{\art}[6]{{\sc #1, \rm #2, \it #3 \bf #4 \rm (#5), \mbox{#6}.}}
\newcommand{\auth}[2]{{#1, #2.}}

\newcommand{\artprep}[3]{{\sc #1, \rm #2, #3.}}
\newcommand{\arttoappear}[3]{{\sc #1, \rm #2, to appear in \it #3}}
\newcommand{\book}[3]{{\sc #1, \it #2, \rm #3.}}
\newcommand{\AND}{{\rm and }}
\RequirePackage{amsthm}
\newtheoremstyle{descriptive}%
  {\topsep}   
  {\topsep}   
  {\rmfamily} 
  {}          
  {\bfseries} 
  {.}         
  { }         
  {}          
\newtheoremstyle{propositional}%
  {\topsep}   
  {\topsep}   
  {\itshape}  
  {}          
  {\bfseries} 
  {.}         
  { }         
  {}          
\theoremstyle{propositional}
\newtheorem{thm}{Theorem}[section]
\newtheorem{prop}[thm]{Proposition}
\newtheorem{lem}[thm]{Lemma}
\newtheorem{cor}[thm]{Corollary}
\theoremstyle{descriptive}
\newtheorem{deff}[thm]{Definition}
\newtheorem{example}[thm]{Example}
\newtheorem{remark}[thm]{Remark}
\makeatletter
\renewenvironment{proof}[1][\proofname]{\par
  \pushQED{\qed}%
  \normalfont
  \trivlist
  \item[\hskip\labelsep
        \itshape
    #1\@addpunct{.}]\ignorespaces
}{%
  \popQED\endtrivlist\@endpefalse
} \makeatother
\newcommand{\setm}{\setminus}
\renewcommand{\emptyset}{\varnothing}
\newcommand{\Cp}{{C_p}}
\DeclareMathOperator{\capp}{cap}
\newcommand{\cp}{\capp_p}
\DeclareMathOperator{\diam}{diam} 
\DeclareMathOperator{\Lip}{Lip}
\newcommand{\Lipc}{{\Lip_c}}
\DeclareMathOperator{\unif}{unif}

\DeclareMathOperator*{\osc}{osc}
\DeclareMathOperator*{\dvg}{div}
\DeclareMathOperator*{\essinf}{ess\,inf}
\DeclareMathOperator*{\esssup}{ess\,sup}
\DeclareMathOperator*{\finelimsup}{fine\,lim\,sup}
\DeclareMathOperator*{\fineliminf}{fine\,lim\,inf}
\DeclareMathOperator*{\finelim}{fine\,lim}
\DeclareMathOperator*{\esslim}{ess\,lim}
\DeclareMathOperator{\fineint}{fine-int}
\newcommand{\cpessinf}{\text{$\Cp$-}\essinf}
\DeclareMathOperator*{\cpessinfalt}{\text{$\Cp$-}\essinf}
\newcommand{\cplim}{\text{$\Cp$-}\esslim}
\DeclareMathOperator*{\cplimalt}{\text{$\Cp$-}\esslim}

\newcommand{\pSubset}{\mathrel{\overset{p}{\Subset}}}

\newcommand{\bdry}{\partial}
\newcommand{\bdy}{\bdry}
\newcommand{\loc}{_{\rm loc}}
{\catcode`p =12 \catcode`t =12 \gdef\eeaa#1pt{#1}}      
\def\accentadjtext#1{\setbox0\hbox{$#1$}\kern   
                \expandafter\eeaa\the\fontdimen1\textfont1 \ht0 }
\def\accentadjscript#1{\setbox0\hbox{$#1$}\kern 
                \expandafter\eeaa\the\fontdimen1\scriptfont1 \ht0 }
\def\accentadjscriptscript#1{\setbox0\hbox{$#1$}\kern   
                \expandafter\eeaa\the\fontdimen1\scriptscriptfont1 \ht0 }
\def\accentadjtextback#1{\setbox0\hbox{$#1$}\kern       
                -\expandafter\eeaa\the\fontdimen1\textfont1 \ht0 }
\def\accentadjscriptback#1{\setbox0\hbox{$#1$}\kern     
                -\expandafter\eeaa\the\fontdimen1\scriptfont1 \ht0 }
\def\accentadjscriptscriptback#1{\setbox0\hbox{$#1$}\kern 
                -\expandafter\eeaa\the\fontdimen1\scriptscriptfont1 \ht0 }
\def\itoverline#1{{\mathsurround0pt\mathchoice
        {\rlap{$\accentadjtext{\displaystyle #1}
                \accentadjtext{\vrule height1.593pt}
                \overline{\phantom{\displaystyle #1}
                \accentadjtextback{\displaystyle #1}}$}{#1}}
        {\rlap{$\accentadjtext{\textstyle #1}
                \accentadjtext{\vrule height1.593pt}
                \overline{\phantom{\textstyle #1}
                \accentadjtextback{\textstyle #1}}$}{#1}}
        {\rlap{$\accentadjscript{\scriptstyle #1}
                \accentadjscript{\vrule height1.593pt}
                \overline{\phantom{\scriptstyle #1}
                \accentadjscriptback{\scriptstyle #1}}$}{#1}}
        {\rlap{$\accentadjscriptscript{\scriptscriptstyle #1}
                \accentadjscriptscript{\vrule height1.593pt}
                \overline{\phantom{\scriptscriptstyle #1}
                \accentadjscriptscriptback{\scriptscriptstyle #1}}$}{#1}}}}
\def\itunderline#1{{\mathsurround0pt\mathchoice
        {\rlap{$\underline{\phantom{\displaystyle #1}
                \accentadjtextback{\displaystyle #1}}$}{#1}}
        {\rlap{$\underline{\phantom{\textstyle #1}
                \accentadjtextback{\textstyle #1}}$}{#1}}
        {\rlap{$\underline{\phantom{\scriptstyle #1}
                \accentadjscriptback{\scriptstyle #1}}$}{#1}}
        {\rlap{$\underline{\phantom{\scriptscriptstyle #1}
                \accentadjscriptscriptback{\scriptscriptstyle #1}}$}{#1}}}}
\newcommand{\de}{\delta}
\newcommand{\la}{\lambda}
\newcommand{\ga}{\gamma}
\newcommand{\Ga}{\Gamma}
\newcommand{\Om}{\Omega}
\renewcommand{\phi}{\varphi}
\newcommand{\p}{{$p\mspace{1mu}$}}
\newcommand{\clEp}{{\itoverline{E}\mspace{1mu}}^p}
\newcommand{\clVp}{{\overline{V}\mspace{1mu}}^p}
\newcommand{\bdyp}{\bdy_p} 
\newcommand{\R}{\mathbf{R}}
\newcommand{\eR}{{\overline{\R}}}
\newcommand{\Q}{\mathbf{Q}}
\newcommand{\K}{{\cal K}}
\newcommand{\Np}{N^{1,p}}
\newcommand{\Lp}{L^p}
\newcommand{\Npploc}{N^{1,p}_{\textup{fine-loc}}}
\newcommand{\Lploc}{L^{p}\loc}
\newcommand{\grad}{\nabla}
\newcommand{\Et}{\widetilde{E}}
\newcommand{\sS}{{S}}
\newcommand{\sQ}{{Q}}
\newcommand{\sP}{{P}}
\newcommand{\sR}{{R}}
\newcommand{\uQ}{{\itoverline{Q}}}
\newcommand{\lQ}{{\itunderline{Q}}}
\newcommand{\uS}{{\itoverline{S}}}
\newcommand{\uR}{{\itoverline{R}}}
\newcommand{\lS}{{\itunderline{S}}}
\newcommand{\lR}{{\itunderline{R}}}
\newcommand{\uP}{{\itoverline{P}}}
\newcommand{\lP}{{\itunderline{P}}}
\newcommand{\UU}{{\cal U}}
\newcommand{\UUt}{\widetilde{\cal U}}
\newcommand{\uShat}{{\widehat{S}}}
\newcommand{\uPhat}{{\widehat{P}}}
\newcommand{\uQhat}{{\widehat{Q}}}
\newcommand{\uRhat}{{\widehat{R}}}
\newcommand{\ut}{\tilde{u}}
\newcommand{\vt}{\tilde{v}}
\newcommand{\reg}{\sharp}
\newcommand{\Cunif}{C_{\unif}}
\newcommand{\tP}{\widetilde{P}}
\newcommand{\tS}{\widetilde{S}}
\newcommand{\imp}{\ensuremath{\Rightarrow} }
\numberwithin{equation}{section}
\newenvironment{ack}{\medskip{\it Acknowledgement.}}{}

\begin{document}

\authortitle{Anders Bj\"orn, Jana Bj\"orn and Visa Latvala}
            {The Perron method
associated with finely \p-harmonic functions on finely open sets}
\author{
Anders Bj\"orn \\
\it\small Department of Mathematics, Link\"oping University, SE-581 83 Link\"oping, Sweden\\
\it \small anders.bjorn@liu.se, ORCID\/\textup{:} 0000-0002-9677-8321
\\
\\
Jana Bj\"orn \\
\it\small Department of Mathematics, Link\"oping University, SE-581 83 Link\"oping, Sweden\\
\it \small jana.bjorn@liu.se, ORCID\/\textup{:} 0000-0002-1238-6751
\\
\\
Visa Latvala \\
\it\small Department of Physics and Mathematics,
University of Eastern Finland, \\
\it\small  P.O. Box 111, FI-80101 Joensuu,
Finland\/{\rm ;} \\
\it \small visa.latvala@uef.fi, ORCID\/\textup{:} 0000-0001-9275-7331
}

\date{}

\maketitle

\noindent{\small
{\bf Abstract}.}
Given a bounded finely open set $V$ and a function $f$
on the fine boundary of $V$,
we introduce four types of upper Perron solutions
to the nonlinear Dirichlet problem for \p-energy minimizers, $1<p<\infty$,
with $f$ as boundary data.
These solutions are
given as pointwise infima of suitable families of fine \p-superminimizers
in $V$.
We show  (under natural assumptions) that  the four upper
Perron solutions are equal quasieverywhere and that they are fine \p-minimizers
of the \p-energy integral.
We moreover show that the upper and lower Perron solutions coincide
quasieverywhere 
for Sobolev and for uniformly
continuous boundary data, i.e.\ that such boundary data
are resolutive.
For the uniformly continuous boundary data, the Perron
solutions are also shown to be finely continuous and thus finely \p-harmonic.
We prove our results in 
a complete metric space $X$ equipped with a doubling
measure supporting a \p-Poincar\'e inequality, but they are new also in unweighted $\R^n$.

\medskip


\medskip


\section{Introduction}

The aim in this paper is to study the  Perron method related to the
nonlinear Dirichlet problem for \emph{finely} \p-harmonic functions  on
\emph{finely open} sets.
Recall that in unweighted $\R^n$, \p-harmonic functions
on open sets
minimize the \p-energy $\int |\grad u|^p\,dx$ and are
solutions of the \p-Laplace equation
\begin{equation}  \label{eq-p-Lapl}
\Delta_p u := \dvg(|\grad u|^{p-2} \grad u) =0.
\end{equation}
The \p-superharmonic functions on open sets
associated with this equation are known to
be finely continuous, which has (especially in the 
classical linear case $p=2$)
led to the development of fine potential theory.
For $p\ne2$, the equation~\eqref{eq-p-Lapl} is nonlinear.

Given a bounded finely open set $V$ (with 
complement of positive capacity)
and a function $f$
on the fine boundary $\bdyp V$,
we seek a finely
\p-harmonic function in $V$ with $f$ as boundary data.
For this purpose, we introduce four types of upper Perron solutions
$\uP f$, $\uQ f$, $\uR f$ and $\uS f$,
given as pointwise infima of suitable families of fine superminimizers in $V$.
Out of these, $\uQ f$ is the most important for our approach as it,
surprisingly, combines several useful properties of the solutions.
On the other hand, $\uP f$ most closely mimics the traditional Perron solutions.
The Perron solutions $\uS f$ and $\uR f$ provide useful connections
between the other two solutions.

We perform our study in a complete metric space $X$ equipped with a doubling
measure supporting a \p-Poincar\'e inequality, where $1<p<\infty$.
However, as far as we know,
the Perron method associated with \p-harmonic functions has not
been studied earlier on finely open sets
beyond the linear case $p=2$, not even on unweighted $\R^n$.
The fine (super)minimizers considered in this paper
coincide on unweighted $\R^n$ with the fine \p-(super)solutions
of $\Delta_p u =0$
introduced in Kilpel\"ainen--Mal\'y~\cite{KiMa92} in 1992.

The following is a special case of our main results and follows directly from
Theorems~\ref{thm-all-equal} and~\ref{thm-uQ-finemin}.
  By "q.e." we mean "quasieverywhere", i.e.\ up to a set of zero capacity.

\begin{thm} \label{thm-main-Perron}
Let $f:\bdy_p V \to \R$ be bounded. Then
\[
\uP f=\uQ f =\uR f=\uS f
\quad \text{q.e.\ in $V$}
\]
and
the four upper Perron solutions are fine minimizers in $V$.
\end{thm}

The minimizing property of the Perron solutions is vital but nontrivial.
Both parts of Theorem~\ref{thm-main-Perron}
hold also for large classes of unbounded functions,
see Theorems~\ref{thm-all-equal} and~\ref{thm-uQ-finemin}
for precise statements.

On open sets, Perron solutions are defined using
\p-superharmonic functions, which essentially are lsc-regularized superminimizers
of the \p-energy integral.
Here we have defined the upper Perron solutions $\uP f$ and $\uS f$ using
finely lsc-regularized fine superminimizers.
In contrast, the upper Perron solutions $\uQ f$ and $\uR f$ are  defined
using finely usc-regularized fine superminimizers.
Note that (unlike in the linear case $p=2$ or on open sets) 
it is not known whether fine (super)minimizers have 
finely continuous representatives, not even in unweighted $\R^n$.

It is perhaps a bit surprising that we have obtained the best results using
usc-regularizations, but
we have not been able to construct direct proofs showing e.g.\ that
$\uP f$ is a fine minimizer or that $\lP f\le \uP f$ q.e.
On the other hand, $\uQ f$ and $\uR f$ are finely upper semicontinuous and can be
written q.e.\ as decreasing limits of fine superminimizers 
(Lemma~\ref{lem-uQ-defining-seq}).
This is then used as a key tool to obtain several properties
of the Perron solutions, as well as comparisons between them,
  including Theorem~\ref{thm-main-Perron}.
Our proofs are therefore rather different from the corresponding proofs for open sets
in the literature, where
the compactness of the boundary or
convergence on countable dense subsets is used.
Since countable sets often have zero capacity  and are not seen 
by the fine topology, they cannot be used in our situation.

In order to use the Perron method, it is also important
to study resolutivity, i.e.\
when the upper and lower Perron solutions agree (at least q.e.).
The following result is a special case of
Proposition~\ref{prop-h-Qf-new} and
Theorems~\ref{thm-all-equal} and~\ref{thm-cont-Perron}.
In what follows, $\clVp = V \cup \bdyp V$ is the fine closure of $V$
and $\Np(V)$ is the Newtonian--Sobolev space, defined by means 
of upper gradients, see Section~\ref{sect-prelim}.

\begin{thm}  \label{thm-h-Qf-intro}
Let $f : \clVp \to [-\infty,\infty]$.
\begin{enumerate}
\item
If $f\in\Np(V)$ and
\[
  f(x)=\finelim_{V\ni y\to x}f(y)
  \quad \text{for q.e.\ $x\in\bdy_p V$}
\]
\textup{(}which in particular holds if $f\in\Np(X)$\textup{),}
then
\begin{equation} \label{eq-h-intro}
  \lP f = \uP f = \lQ f = \uQ f=\lR f = \uR f=\lS f = \uS f 
  \quad \text{q.e.\ in } V.
\end{equation}
\item
If $f$  is uniformly continuous
on $\bdyp V$
  with respect to the metric topology, then
  $\lQ f = \uQ f =\lS f = \uS f $ everywhere in $V$,
  and this solution is finely continuous in $V$.
\end{enumerate}
\end{thm}

If we ignore the ``q.e.'' in \eqref{eq-h-intro},
Theorem~\ref{thm-h-Qf-intro} corresponds to the
  two main resolutivity results known in the nonlinear case on \emph{open sets}:
  for functions in $\Np$ and for continuous functions.
There are also a few other resolutivity results on open sets, mainly
for semicontinuous functions. (See e.g.\ \cite[Chapter~10]{BBbook}, or
the references below.)
When turning to finely open sets
there are
new issues involved and many open questions,
see e.g.\ Section~9 in our paper~\cite{BBLat4}.

Another question is invariance of the Perron solutions under
perturbations of the boundary data on sets of zero capacity.
On open sets, such invariance results for nonlinear problems
were first obtained
by Bj\"orn--Bj\"orn--Shan\-mu\-ga\-lin\-gam~\cite{BBS2} for
functions in
$\Np$ and for continuous functions.
In our situation, such invariance for the $\sQ$- and $S$-Perron solutions follows
immediately from their definitions.
Moreover,
these two
upper Perron solutions can equivalently be defined
using (finely usc/lsc-regularized)  \emph{fine minimizers}
(see Remark~\ref{rmk-h}).

Perron solutions of the Dirichlet problem were introduced by
Perron~\cite{Perron23} and Remak~\cite{remak} (independently) in 1923
for harmonic functions on open subsets of the plane.
The Perron method
has later been 
used also for  
various linear and nonlinear elliptic and parabolic equations.
Fuglede~\cite[p.~173]{Fug} extended the Perron method
to finely open sets in the linear axiomatic setting.
An alternative approach, under weaker axioms,
was
   given by Luke\v{s}--Mal\'y--Zaj\'i\v{c}ek~\cite[Chapter~13]{LuMaZa}.
They however needed to consider boundary values on the ``bitopological'' boundary
$\bdy V \setm V$.
Our approach is thus more in line with  Fuglede's.

However, the linear theory
considered in~\cite{Fug} and~\cite{LuMaZa}
differs much from the nonlinear one since it is based on potentials
and measures rather than on Sobolev spaces.
At the same time for nonlinear problems, such as~\eqref{eq-p-Lapl},
Sobolev spaces are
essential already for defining \p-harmonic functions
on $\R^n$, as seen in
Heinonen--Kilpel\"ainen--Martio~\cite{HeKiMa}.
The metric space approach is particularly suitable
  for nonlinear fine potential theory
since Newtonian (Sobolev) spaces can be defined  directly
on an arbitrary measurable subset $E$ by  considering
$E$ as a metric space in its own right.
On the other hand, it
follows from Bj\"orn--Bj\"orn~\cite[Theorem~7.3]{BBnonopen}
and Bj\"orn--Bj\"orn--Latvala~\cite{BBLat4}
that  the study of
the Dirichlet problem for \p-harmonic functions
is  not natural beyond
finely open or quasiopen sets.
By Theorem~\ref{thm-finelyopen-quasiopen},
finely open sets and superlevel
sets of global Newtonian functions are quasiopen.
The use of fine limits in our definition of Perron solutions
makes it natural to restrict attention to
finely open sets.

As far as we know, Aronsson~\cite[Section~4]{aronsson67} was
the first who used the Perron method in connection with a nonlinear
equation, namely for the $\infty$-Laplacian (also called
Aronsson's equation).
In connection with \p-harmonic functions (for $1<p<\infty$)
it was first used by
Granlund--Lindqvist--Martio~\cite{GLM86}.
See also Kilpel\"ainen~\cite{Kilp89},
Heinonen--Kilpel\"ainen--Martio~\cite{HeKiMa},
Bj\"orn--Bj\"orn--Shan\-mu\-ga\-lin\-gam~\cite{BBS2}, \cite{BBSdir},
Bj\"orn--Bj\"orn--Sj\"odin~\cite{BBSjodin} and
Hansevi~\cite{hansevi2}.
All of these studies were on open sets.
The nonlinear Dirichlet problem for fine
  solutions/minimizers
on quasiopen sets was studied in the Sobolev sense by
Kilpel\"ainen--Mal\'y~\cite{KiMa92} on unweighted $\R^n$,
and on metric spaces by the authors
in~\cite{BBLat4}.

The outline of the paper is as follows.
In Section~\ref{sect-prelim} we present the main background definitions
related to Newtonian functions and capacity in metric spaces, while in
Section~\ref{sect-fine-cont} we introduce the fine topology
and finely semicontinuous regularizations.
Following our   papers~\cite{BBLat3} and~\cite{BBLat4},
we introduce the fine local Newtonian space $\Npploc$, fine superminimizers and the fine obstacle problem
in Section~\ref{sect-fine-min}.
The main new contribution here is Proposition~\ref{prop-3.2}, which
is a convergence result for solutions of obstacle problems  that is important
later on.

In Section~\ref{sect-fine-perron} we introduce our Perron solutions and prove
two auxiliary results for the upper Perron solution $\uQ f$,
whereas Section~\ref{sect-comparison} contains comparisons
between different
Perron solutions.
In Section~\ref{sect-Perron-fine-min} we give sufficient conditions
for Perron solutions to be fine minimizers.
Finally in Section~\ref{sect-continuity} we
prove Theorem~\ref{thm-h-Qf-intro}
and in Section~\ref{sect-open}
we study 
the behaviour of our Perron solutions on open sets.

\begin{ack}
A.~B. and J.~B. were supported by the Swedish Research Council,
grants 2016-03424 and 2020-04011 resp.\ 621-2014-3974 and 2018-04106.
Part of this research was done during three visits
of V.~L. to Link\"oping University in
2017--2022
and a visit, supported by SVeFUM, of A.~B. and J.~B. to the University
of Eastern Finland in Joensuu in 2022.
We thank all these institutions for their hospitality and support.
\end{ack}

\section{Notation and preliminaries}
\label{sect-prelim}

In this section, we introduce
the necessary metric space concepts used in this paper.
For brevity, we refer to
our papers~\cite{BBLat1} and~\cite{BBLat2} for more
extensive introductions 
and references to the literature.
See also the monographs Bj\"orn--Bj\"orn~\cite{BBbook} and
Heinonen--Koskela--Shan\-mu\-ga\-lin\-gam--Tyson~\cite{HKST},
where the theory of upper gradients and Newtonian (Sobolev) spaces
on metric spaces is thoroughly  developed with proofs.

Let $X$ be a metric space equipped
with a metric $d$ and a positive complete  Borel  measure $\mu$
such that $\mu(B)<\infty$ for all balls $B \subset X$.
We also assume that $1<p< \infty$.

A \emph{curve} is a continuous mapping from an interval,
and a \emph{rectifiable} curve is a curve with finite length.
We will only consider curves which are nonconstant, compact and rectifiable,
and they can therefore be parameterized by their arc length $ds$.
A property holds for \emph{\p-almost every curve}
if the curve family $\Ga$ for which it fails has zero \p-modulus,
i.e.\ there is $\rho\in L^p(X)$ such that
$\int_\ga \rho\,ds=\infty$ for every $\ga\in\Ga$.

  A measurable
  function $g:X \to [0,\infty]$ is a \p-weak \emph{upper gradient}
of $u:X \to \eR:=[-\infty,\infty]$
if for \p-almost all curves
$\gamma: [0,l_{\gamma}] \to X$,
\begin{equation*} 
        |u(\gamma(0)) - u(\gamma(l_{\gamma}))| \le \int_{\gamma} g\,ds,
\end{equation*}
where the left-hand side is $\infty$
whenever at least one of the
terms therein is infinite.
If $u$ has a \p-weak upper gradient in $\Lploc(X)$, then
it has a \emph{minimal \p-weak upper gradient}
$g_u \in \Lploc(X)$
in the sense that
$g_u \le g$ a.e.\
for every \p-weak upper gradient $g \in \Lploc(X)$ of $u$.
For such measurable $u$, we let
\[
        \|u\|_{\Np(X)} = \biggl( \int_X |u|^p \, d\mu
                +  \int_X g_u^p \, d\mu \biggr)^{1/p}.
\]
The \emph{Newtonian space} on $X$ is
\[
        \Np (X) = \{u: \|u\|_{\Np(X)} <\infty \}.
\]
For a measurable set $E\subset X$, the Newtonian space $\Np(E)$
is
defined by
considering $(E,d|_E,\mu|_E)$ as a metric space in its own right.

The space $\Np(X)/{\sim}$, where  $u \sim v$ if and only if $\|u-v\|_{\Np(X)}=0$,
is a Banach space and a lattice.
In this paper it will be convenient to
assume that functions in $\Np(X)$
 are defined everywhere (with values in $\eR$),
not just up to an equivalence class in the corresponding function space.
For an arbitrary set $A \subset X$,  we let
\[
  \Np_0(A)=\{u|_{A} : u \in \Np(X) \text{ and }
  u=0 \text{ on } X \setm A\}.
\]
Functions from $\Np_0(A)$ can be extended by zero in $X\setm A$ and we
will regard them in that sense when
needed.

The  \emph{Sobolev capacity} of an arbitrary set $A\subset X$ is
\[
 \Cp(A)=\inf_{u}\|u\|_{\Np(X)}^p,
\]
where the infimum is taken over all $u \in \Np(X)$ such that
$u\geq 1$ on $A$.
A property holds \emph{quasieverywhere} (q.e.)\
if the set of points  for which it fails has capacity zero.
The capacity is the correct gauge
for distinguishing between two Newtonian functions,
  namely $\|u\|_{\Np(X)}=0$ if and only if $u=0$ q.e.
Moreover, if $u,v \in \Np(X)$ and $v= u$ a.e., then $v=u$ q.e.

The measure  $\mu$  is \emph{doubling} if
there is $C>0$ such that for all balls
$B=B(x_0,r):=\{x\in X: d(x,x_0)<r\}$ in~$X$,
we have
$        0 < \mu(2B) \le C \mu(B) < \infty$,
where $\lambda B=B(x_0,\lambda r)$.
In this paper, all balls are open.

The space $X$ supports a \emph{\p-Poincar\'e inequality} if
there are $C>0$ and $\lambda \ge 1$
such that for all balls $B \subset X$,
all integrable functions $u$ on $X$ and all \p-weak upper gradients $g$ of $u$,
we have
$0<\mu(B)<\infty$ and
\begin{equation} \label{PI-ineq}
       \frac{1}{\mu(B)} \int_{B} |u-u_B| \,d\mu
       \le C \diam(B) \biggl( \frac{1}{\mu(\la B)}
              \int_{\lambda B} g^{p} \,d\mu \biggr)^{1/p},
\end{equation}
where $ u_B := \int_B u\, d\mu/\mu(B)$.

In $\R^n$ equipped with a doubling measure $d\mu=w\,dx$,
the \p-Poincar\'e inequality~\eqref{PI-ineq}
is equivalent to the \emph{\p-admissibility} of the weight $w$ in the
sense of Heinonen--Kilpel\"ainen--Martio~\cite{HeKiMa}, see
Corollary~20.9 in~\cite{HeKiMa}
and Proposition~A.17 in~\cite{BBbook}.
Moreover, in this case $g_u=|\nabla u|$ a.e.\ if $u \in \Np(\R^n,\mu)$,
and the capacities in this paper
coincide
with the corresponding capacities in~\cite{HeKiMa}
(see \cite[Theorem~6.7]{BBbook} and \cite[Theorem~5.1]{BBvarcap}).

\section{Fine topology}
\label{sect-fine-cont}

\emph{Throughout the rest of the paper, we assume
  that $X$ is complete and  supports a \p-Poincar\'e inequality
with  $1<p< \infty$ and
that
$\mu$ is doubling.
We also assume that $V$ is a nonempty finely open set,
and from Section~\ref{sect-fine-min} onwards
that in addition $V$ is bounded and
$\Cp(X \setm V)>0$.}

\medskip

To avoid pathological situations we also assume that $X$
contains at least two points (and thus must be uncountable due to the
Poincar\'e inequality).
In this section we recall the basic facts about the fine topology
associated with Newtonian functions.
The \emph{variational capacity} 
of $A$ with respect to $B$ is defined by
\[
\cp(A,B) := \inf_{u}\int_{X} g_{u}^p\, d\mu,
\]
where the infimum is taken over all $u \in \Np_0(B)$
such that $u\geq 1$ on $A$.

\begin{deff}\label{deff-thinness}
A set $E\subset X$ is  \emph{thin} at $x\in X$ if
\begin{equation*} 
\int_0^1\biggl(\frac{\cp(E\cap B(x,r),B(x,2r))}{\cp(B(x,r),B(x,2r))}\biggr)^{1/(p-1)}
     \frac{dr}{r}<\infty.
\end{equation*}
A set $V\subset X$ is \emph{finely open} if
$X\setminus V$ is thin at each point $x\in V$.
\end{deff}

In the definition of thinness,
we use the convention that the integrand
is 1 whenever $\cp(B(x,r),B(x,2r))=0$.
It is easy to see that the finely open sets give rise to a
topology, which is called the \emph{fine topology}.
Every open set is finely open, but the converse is not true in general.
A function $u : V \to \eR$, defined on a finely open set $V$, is
\emph{finely continuous} if it is continuous when $V$ is equipped with the
fine topology and $\eR$ with the usual topology.
Pointwise fine (semi)continuity is defined analogously.
The
fine interior, fine boundary and fine closure of $E$
are denoted $\fineint E$, $\bdyp E$ and $\clEp$, respectively.
See Bj\"orn--Bj\"orn~\cite[Section~11.6]{BBbook}
and Bj\"orn--Bj\"orn--Latvala~\cite{BBLat1}
for further discussion on thinness and the fine topology
in metric spaces. Note that the fine topology is usually
not metrizable.

A set $U\subset X$ is \emph{quasiopen} if for every
$\varepsilon>0$ there is an open set $G\subset X$ such that $\Cp(G)<\varepsilon$
and $G\cup U$ is open.
Quasiopen sets are measurable by Lemma~9.3 in
  Bj\"orn--Bj\"orn~\cite{BBnonopen}.
Various characterizations of  quasiopen sets can be found in
Bj\"orn--Bj\"orn--Mal\'y~\cite{BBMaly}.
Therein it is also shown that the $\Cp$ capacities with respect to $X$ and
 a quasiopen $U$ have the same zero sets.

The following result explains the close connection between
finely open and quasiopen sets.
In particular, 
finely open sets are quasiopen and we therefore have at our disposal
all earlier results obtained for quasiopen sets.

\begin{thm} \label{thm-finelyopen-quasiopen}
  \textup{(Theorem~3.4 in~\cite{BBLat4})} 
The following conditions are equivalent for any set $U\subset X$\/\textup{:}
\begin{enumerate}
\item \label{tt-a}
$U$ is quasiopen\textup{;}
\item 
$U=V \cup E$ for some finely open $V$ and a set $E$
with $\Cp(E)=0$\textup{;}
\item \label{t-d}
$\Cp(U\setm \fineint U)=0$\textup{;}
\item \label{t-dd}
$U=\{x:u(x)>0\}$ for some $u\in \Np(X)$.
\end{enumerate}
\end{thm}

We will use the following principle several times 
as a substitute for compactness.

\begin{thm} \label{thm-quasiLindelof}
\textup{(Quasi-Lindel\"of principle, Theorem~3.4 in~\cite{BBLat3})}
 For each family $\mathcal V$ of finely open sets there is a countable subfamily $\mathcal V'$ such that
\[
\Cp\biggl(\bigcup_{V\in\mathcal V}V\setminus \bigcup_{V'\in\mathcal V'}V'\biggr)=0.
\]
\end{thm}

We will also need the following characterization of $\Np_0(V)$, which
is a special case of Theorem~7.2 in~\cite{BBLat4}.

\begin{prop}   \label{prop-Np0}
Let $u \in \Np(V)$.
Then $u \in \Np_0(V)$ if and only if
\[
  \finelim_{V \ni y\to x} u(y)=0
\quad \text{for q.e.\ } x\in \bdyp V.
\]
\end{prop}

In this paper, the notions of $\finelim$, $\finelimsup$ and
$\fineliminf$ are defined using punctured fine neighbourhoods.
Since
\[
\cp(B(x,r) \setm\{x\},B(x,2r))= \cp(B(x,r),B(x,2r)),
\]
there are no isolated points in the fine topology,
i.e.\ no singleton sets are finely open.
Moreover,  if $W$ is finely open and $\Cp(E)=0$,
then $W \setm E$ is also finely open.
Hence,
if $x \in \clVp$, then 
\[
\fineliminf_{V \ni y \to x} u(y)
  = \sup_{\text{finely open } W\ni x}  \, \inf_{(V\cap W) \setm \{x\}} u
 = \sup_{\text{finely open } W\ni x}  \, \cpessinfalt_{(V\cap W) \setm \{x\}} u,
\]
where
\[
 \cpessinf_{E} u \\
 :=  \sup \{ k: u \ge k \text{ q.e.\ in } E \}.
 \]

We will
 extensively use \emph{fine lsc- and usc-regularizations}, which we define as follows.
Let $u:V \to \eR$ and set for $x \in V$,
\begin{align*}
u_\reg(x) &= \fineliminf_{V \ni y \to x} u(y)
 = \sup_{\text{finely open } W\ni x} \, \cpessinfalt_{(V\cap W) \setm \{x\}} u,
 \nonumber\\
u_*(x) &= \begin{cases}
  u_\reg(x), & \text{if } \Cp(\{x\})=0, \\
  \min\{u_\reg(x),u(x)\}, & \text{if } \Cp(\{x\})>0,
\end{cases}
\quad
 = \sup_{\text{finely open } W\ni x} \,  \cpessinf_{V\cap W} u.
\end{align*}
Thus, $u_\reg$ is defined using punctured neighbourhoods, while
$u_*$ is using nonpunctured neighbourhoods.
We also say that $u$ is \emph{finely lsc-regularized} if $u=u_*$.
Similarly, we define
\begin{align*}
  u^\reg(x) &= \finelimsup_{V \ni y \to x} u(y)
\quad   \text{and} \quad
u^*(x) = \begin{cases}
  u^\reg(x), & \text{if } \Cp(\{x\})=0, \\
  \max\{u^\reg(x),u(x)\}, & \text{if } \Cp(\{x\})>0,
\end{cases}
\end{align*}
and say that $u$ is \emph{finely usc-regularized} if $u=u^*$.

It follows directly from the definition that
for every $a\in\R$ the set
$\{x\in V: u_*(x)>a\}$ is finely open and hence $u_*$ is finely lower
semicontinuous in $V$.
Moreover,  $u$ is finely lower semicontinuous in $V$ if and only if
$u \le u_*$ everywhere in $V$.
The same is true for $u_\reg$ and similar results hold
for $u^*$ and  $u^\reg$.

\begin{remark}\label{rmk-reg}
Assume that $u$ and $v$ are functions on $V$
such that $u\le v$ q.e.\ in $V$.
Then
\[
u_*\le v_*,
\
u_\reg\le v_\reg,
\
u^*\le v^*
\text{ and }
u^\reg\le v^\reg
\quad \text{everywhere in }V.
\]
In fact, if $\Cp(\{x\})=0$, then
\[
u_*(x)=u_\reg(x)\le v_\reg(x)=v_*(x)
\]
since the fine lower limits do not see sets of zero capacity.
If $\Cp(\{x\})>0$, then $u(x)\le v(x)$ by assumption,
and $u_\reg(x)\le v_\reg(x)$ as above. The upper semicontinuous case is similar.
\end{remark}

The following lemma is the reason for introducing $u_*$ with the more complicated
definition than $u_\reg$, since $u_{\reg\reg}=u_\reg$ is false in general,
see Example~\ref{ex-u-reg} below.
Note however that  most of the functions
we regularize are finely continuous q.e., and thus $u_\reg=u_*$.
In particular this holds for fine superminimizers and functions in
Newtonian spaces.
It is only when we regularize the Perron solutions that
the notion of $u_*$ is really needed.

\begin{lem}   \label{lem-u*-fine-reg}
  The function $u_*$ is finely lsc-regularized in $V$, i.e.\ $u_{**}=u_*$.
Similarly $u^{**}=u^*$.
\end{lem}

\begin{proof}
It suffices to show that $u_*\ge u_{**}$ in $V$, since the converse inequality
follows from the fine lower semicontinuity of $u_*$. Let $x\in V$ and
\begin{equation}   \label{eq-a<fineliminf-new}
a < u_{**}(x)
\end{equation}
be arbitrary.
By definition, there exists a finely open set $W\ni x$ such that
$u_*>a$ q.e.\ in $W$.
Hence for q.e.\  $y\in W$, there exists a finely open
set $W_y\ni y$ such that $u>a$ q.e.\ in $W_y$.

The sets $W_y$ cover $W$ up to a set of zero capacity.
Therefore, by the quasi-Lindel\"of property
(Theorem~\ref{thm-quasiLindelof}),
there are $\{y_j\}_{j=1}^\infty$ such that the sets $W_{y_j}$ cover $W$ up to
a set of zero capacity.
Since $\{y_j\}_{j=1}^\infty$ is countable, we conclude that $u>a$
q.e.\ in $W$.
From this it follows that $u_*(x)\ge a$ and taking supremum over all $a$
admissible in \eqref{eq-a<fineliminf-new} gives $u_*\ge u_{**}$.
The proof for $u^*$ is similar.
\end{proof}

\begin{example} \label{ex-u-reg}
We shall now see that the statement corresponding to Lemma~\ref{lem-u*-fine-reg}
is false in general for $u_\reg$.
Consider
  the
real line $\R$.
Since $p>n=1$, all points have positive capacity and the fine topology
is the same as the Euclidean topology.
Let
$u(2^{-n})=0$, $n=1,2,\ldots$\,,
and $u=1$ otherwise.
Then $u_{\reg}(0)=0$ and $u_{\reg}=1$ otherwise,
while $u_{\reg\reg} \equiv 1$, and thus $u_{\reg\reg} \ne u_\reg$.

Similarly, if we let $v(2^{-n} +2^{-n-m})=0$, $n,m=1,2,\ldots$\,,
and $v=1$ otherwise, then $v_{\reg}=u$ and thus
$v_{\reg\reg}=u_{\reg} \not \equiv 1$, while $v_{\reg\reg\reg} = u_{\reg\reg} \equiv 1$.
One can create similar examples such that any fixed number
of $\reg$-lsc-regularizations is not $\reg$-lsc-regularized.
\end{example}

\begin{remark} \label{rmk-u_*-le-u^*}
A direct consequence of the definitions is that
$u_* \le u^*$ (and   $u_\sharp \le u^\sharp$), a fact that we will use extensively
when comparing the different Perron solutions later in the paper.

If $U$ is quasiopen, then points in
$U \setm {\overline{\fineint U}\mspace{1mu}}^p$
  are isolated in the 
  fine topology,
  and it is thus not natural to define fine
limits and fine regularizations at
such points.
As a simple example,
let
$U \subset \R^n$ (unweighted)
be the union of an open set $G$ and a point
  $x \notin \itoverline{G}$,
  where $1 < p \le n$.

This is another reason for why we have restricted our attention
to finely open sets in this paper.
\end{remark}

\begin{prop} \label{prop-cp-bdyp-U}
Assume that $\Cp(\la B\cap \bdyp V)=0$ for some ball
$B$ with $B \cap V \ne \emptyset$, where
$\la$ is the dilation constant
in the \p-Poincar\'e inequality.
Then $\Cp(B \setm V)=0$.

In particular,
$\Cp(\bdyp V)>0$ if and only if $\Cp(X \setm V)>0$.
\end{prop}

We will not use Proposition~\ref{prop-cp-bdyp-U} directly,
but  $\Cp(\bdyp V)>0$ is important for the definitions
of the $\sS$- and $\sQ$-Perron solutions (see Definition~\ref{def-P})
to make sense
and the comparison principle (Theorem~\ref{thm-qe-comparison-refined})
to hold.

\begin{proof}
Since $\Cp(\la B\cap \bdyp V)=0$, both $\la B\setm V$ and  $\la B\cap V$
are quasiopen by Theorem~\ref{thm-finelyopen-quasiopen}.
It then follows from
Shan\-mu\-ga\-lin\-gam~\cite[Remark~3.5]{Sh-harm}
that for \p-almost all curves $\ga:[0,l_\ga]\to\la B$, both
$\ga^{-1}(V)$ and $[0,l_\ga]\setm \ga^{-1}(V)$ are relatively
open subsets of $[0,l_\ga]$.
Since $[0,l_\ga]$ is connected, $\ga^{-1}(V)$ either equals $[0,l_\ga]$ or
  $\emptyset$.
Hence for \p-almost all curves $\ga\subset\la B$, either
$\ga\subset V$ or $\ga\cap V=\emptyset$.

This implies that the characteristic function $\chi_V$
has $0$ as a \p-weak upper gradient in $\la B$ and so $\chi_V\in\Np(\la B)$.
The \p-Poincar\'e inequality then implies that $\chi_V$ is a.e.\
constant in $B$ and thus also q.e.\ constant, see
  Section~\ref{sect-prelim}.
Since $B\cap V$ is finely open and nonempty, we have
$\Cp(B\cap V)>0$ and hence $\Cp(B\setm V)=0$,
which proves the first statement.

One implication in the second statement then follows by
the countable subadditivity of the capacity upon
letting $r\to\infty$ in $B=B(x,r)$, while the
converse implication is trivial.
\end{proof}

\section{Fine (super)minimizers and the obstacle problem}
\label{sect-fine-min}

\emph{Recall the standing assumptions from the beginning of
Section~\ref{sect-fine-cont}.
In particular, 
$V$ is a bounded nonempty finely open set
with $\Cp(X \setm V)>0$.}

\medskip

To define fine (super)minimizers, we first need
an appropriate fine local Sobolev space.
Here \p-strict subsets will play a key role, as a substitute for
relatively compact subsets.
The results in this section
also hold for quasiopen sets, even though we
only formulate them for finely open sets.
Recall that $W\Subset V$ if $\overline{W}$ is a compact subset of $V$.

\begin{deff}\label{def-fineloc}
  A set $W \Subset V$  is a \emph{\p-strict subset} of $V$ if
there is a function $\eta \in \Np_0(V)$ such that $\eta =1$ on
$W$.
We will write $W \pSubset V$.

A function $u$ belongs to $\Npploc(V)$ if $u \in \Np(W)$
for all finely open $W \pSubset V$.
\end{deff}

By Lemma~3.3 in our paper~\cite{BBLat3},
$V$ has a base of fine neighbourhoods $W \pSubset V$.
Functions in $\Npploc(V)$ are
finite q.e., finely continuous q.e.\ and
quasicontinuous, by Theorem~4.4 in~\cite{BBLat3}.
Throughout the paper, we consider minimal
\p-weak upper gradients in $V$.
For a function $u \in \Npploc(V)$ we say that
$g_{u,V}$ is a \emph{minimal \p-weak upper gradient} of $u$ in $V$
if
\begin{equation*} 
  g_{u,V}=g_{u,W} \text{ a.e.\ in $W$} \quad
  \text{for every finely open  $W \pSubset V$},
\end{equation*}
where $g_{u,W}$ is the
minimal \p-weak upper gradient of $u \in \Np(W)$ with respect to
$W$, defined in Section~\ref{sect-prelim}.
See \cite[Lemma~5.2 and Theorem~5.3]{BBLat3} for the
existence, a.e.-uniqueness and minimality of $g_{u,V}$ when
$u\in\Npploc(V)$.
If $u \in \Np(V)$, then this definition agrees
with the definition  of $g_u$ with respect to $V$
in Section~\ref{sect-prelim}.
If moreover $u\in \Np(X)$, then
the minimal \p-weak upper gradients $g_{u,V}$ and $g_u$ with respect to
$V$ and $X$, respectively, coincide a.e.\ in $V$,
see \cite[Corollary~3.7]{BBnonopen} or \cite[Lemma~4.3]{BBLat3}.
For this reason we drop $V$ from the notation and simply write $g_u$
from now on.

If follows from the definition of $\Npploc(V)$ and the
  properties of $\Np$ (see Section~\ref{sect-prelim}) that
functions in $\Npploc(V)$ are defined everywhere in $V$, and that
if $u \in \Npploc(V)$ and $v=u$ q.e., then also $v \in \Npploc(V)$.
Moreover, by the
quasi-Lindel\"of principle
(Theorem~\ref{thm-quasiLindelof})
and the existence of a base of fine neighbourhoods in~$V$,
Corollary~2.21 in~\cite{BBbook} extends to functions
in $\Npploc(V)$.
That is, if $u,v \in \Npploc(V)$ then
\[
    g_u  = g_v \quad \text{a.e.\ on\/ } \{x \in V : u(x)=v(x)\}.
\]

In~\cite[Definition~5.1]{BBLat4} we introduced the following definition.

\begin{deff}\label{def-finesuper}
A function $u \in \Npploc(V)$ is a
\emph{fine minimizer\/ \textup{(}resp.\ fine
superminimizer\/\textup{)}} in $V$ if
\[
\int_{W} g_{u}^p \, d\mu
\le  \int_{W} g_{u+\phi}^p \, d\mu
\]
for every finely open $W \pSubset V$
and for every (resp.\ every nonnegative)
$\phi \in \Np_0(W)$.
Moreover, $u$ is a \emph{fine subminimizer}
if $-u$ is a fine superminimizer.

A \emph{finely \p-harmonic function} is a finely continuous
fine minimizer.
\end{deff}

On unweighted $\R^n$, it follows from Proposition~5.3 in~\cite{BBLat4}
  that these fine (super)minimizers are exactly the fine \p-(super)solutions
  of $\Delta_p u =0$
  introduced in Kilpel\"ainen--Mal\'y~\cite{KiMa92} in 1992.

Note that, unlike for minimizers and \p-harmonic functions on
  open sets,
\emph{it is not known}, even in unweighted $\R^n$, whether
every fine minimizer
can be modified on a set of zero capacity so that it becomes
  finely continuous and thus finely \p-harmonic,
see the discussion in~\cite[Section~9]{BBLat4}.
Of course,
being a function   from $\Npploc(V)$,
every fine (super/sub)minimizer
is finely continuous at q.e.\ point.

It is not 
difficult to see that a function is a
fine minimizer if and only
if it is both a fine subminimizer and a fine superminimizer,
see Lemma~5.4 in~\cite{BBLat4}.
Corollary~5.8 in \cite{BBLat4} shows that the minimum
of two fine superminimizers is also a fine superminimizer.
We will use these facts without further ado.
We refer to our papers~\cite{BBLat3}, \cite{BBLat4} and~\cite{BBLat5} for
further discussion on fine (super)minimizers
and the Newtonian space $\Npploc(V)$.

In order to prove Theorem~\ref{thm-uQ-finemin},
we will need the following special case
of
the results in~\cite[Section~7]{BBLat5}.

\begin{thm} \label{thm-seq-fine-min-mod}
Let  $\{h_j\}_{j=1}^\infty$ be a decreasing sequence of
fine {\rm(}super\/{\rm)}minimizers in $V$.
Let $h(x)=\lim_{j \to \infty} h_j(x)$.
If either
\[
h \in \Npploc(V)
\quad \text{or} \quad
\esssup_{x \in V} |h(x)-u(x)| < \infty,
\]
for some $u\in\Np(V)$,
then $h$ is a fine {\rm(}super\/{\rm)}minimizer in $V$.
\end{thm}

\begin{proof}
When $h \in \Npploc(V)$ this follows from
\cite[Proposition~7.1 and Corollary~7.3]{BBLat5}.
In the second case, there is $M$ such that $|h-u|\le M$ a.e.
Hence, 
the conclusion follows from~\cite[Theorem~7.4\,(d) and
    Corollary~7.5\,(d)]{BBLat5}
with
$f=u$, $f_0=u-M$ and $f_1=u+M$.
\end{proof}

The obstacle problem is a fundamental tool when
studying fine minimizers.
See \cite{BBnonopen} and~\cite{BBLat4}
for earlier studies of the obstacle problem on nonopen sets
in metric spaces.

\begin{deff} \label{deff-obst-E}
Let $f \in \Np(V)$ and $\psi : V \to \eR$.
Then we define
\begin{equation*}
    \K_{\psi,f}(V)=\{v \in \Np(V) : v-f \in \Np_0(V)
            \text{ and } v \ge \psi \ \text{q.e.\ in } V\}.
\end{equation*}
A function $u \in \K_{\psi,f}(V)$
is a \emph{solution of the $\K_{\psi,f}(V)$-obstacle problem}
if
\begin{equation*} 
       \int_V g^p_{u} \, d\mu
       \le \int_V g^p_{v} \, d\mu
       \quad \text{for all } v \in \K_{\psi,f}(V).
\end{equation*}
\end{deff}

The \emph{Dirichlet problem} is a special case of the obstacle
problem, with the trivial obstacle $\psi \equiv -\infty$.
Note that the boundary data $f$ are only required to belong to $\Np(V)$,
i.e.\ $f$ need not be defined on $\bdry V$
or the fine boundary $\bdyp V$.

\begin{thm}\label{thm-obstacle}
  \textup{(\cite[Theorem~4.2]{BBnonopen} and
    \cite[Theorem~6.2]{BBLat4})}
  Let $f \in \Np(V)$  and $\psi : V \to \eR$,
and assume that $\K_{\psi,f}(V) \ne \emptyset$.
Then there exists a solution $u$ of the $\K_{\psi,f}(V)$-obstacle problem,
which is unique q.e.
The solution $u$ is a fine superminimizer in $V$, and
if $\psi \equiv -\infty$
\textup{(}i.e.\ for the Dirichlet problem\/\textup{)}
  it is a fine minimizer.
\end{thm}

We will also use the following converse of Theorem~\ref{thm-obstacle},
which follows from Lemma~6.3 in
Bj\"orn--Bj\"orn--Latvala~\cite{BBLat5}.

\begin{lem}\label{lem-super-obst}
If $u \in \Np(V)$ is a fine superminimizer in $V$,
then it is a solution of the $\K_{u,u}(W)$-obstacle problem
for any finely open $W\subset V$.
\end{lem}

The following comparison principle was deduced in
Bj\"orn--Bj\"orn~\cite[Corollary~4.3]{BBnonopen} for
more general sets.
We will only need it for finely open  sets.

\begin{lem} \label{lem-obst-le}
\textup{(Comparison principle)}
Let $f,f' \in \Np(V)$
and  $\psi,\psi': V \to \eR$
be such that $\max\{f-f',0\} \in \Np_0(V)$ and
$\psi \le \psi'$ q.e.\ in $V$.
If $u$ and $u'$ are solutions of the
$\K_{\psi,f}(V)$- and  $\K_{\psi',f'}(V)$-obstacle problems,
respectively, then $u \le u'$ q.e.\ in~$V$.
\end{lem}

We end this section with the following
convergence result for solutions of obstacle problems.
It will be a vital tool when proving Lemma~\ref{lem-UUt-UU-seq}.

\begin{prop}\label{prop-3.2}
Let\/ $\{f_j\}_{j=1}^\infty$ be a q.e.\
decreasing sequence in $\Np(V)$ such that
$\|f_j- f_0\|_{\Np(V)}\to 0$, as $j \to\infty$.
Let $u_j$
be a solution of the $\K_{f_j,f_j}(V)$-obstacle problem,
$j=0,1,\ldots$\,.
Then $\{u_j\}_{j=1}^\infty$  decreases q.e.\ to $u_0$.
\end{prop}

This result holds, with the same proof, if $V$
is just assumed to be a bounded measurable set with
$\Cp(X\setm V)>0$.

\begin{proof} As $f_j\ge f_{j+1} \ge f_0$ q.e.\ in $V$, when $j \ge 1$,
it follows from the comparison principle (Lemma~\ref{lem-obst-le}) that
$u_j \ge u_{j+1} \ge u_0$ q.e.\ in $V$.
Hence $\{u_j\}_{j=1}^\infty$ is a q.e.-decreasing sequence
and converges q.e.\ to a function $v$ on $V$.
Moreover, $f_j \to f_0$ q.e., by Corollary~1.72 in~\cite{BBbook}.

Let $w_j=u_j-f_j$ and $w=v-f_0$, all extended by zero outside of $V$.
We have $w_j\in\Np_0(V)\subset\Np(X)$, $j=0,1,\ldots$\,.
Using the
Poincar\'e inequality for $\Np_0$ \cite[Corollary~5.54]{BBbook}
and the fact that $u_j$ is a solution of the $\K_{f_j,f_j}(V)$-obstacle problem,
it follows that
there is $C>0$ such that
\begin{align*}
  \|w_j\|_{\Np(X)}
        &= \bigl( \|w_j\|_{L^p(V)}^p + \|g_{w_j}\|_{L^p(V)}^p\bigr)^{1/p}
       \le C \|g_{w_j}\|_{L^p(V)} \\
       &\le C (\|g_{u_j}\|_{L^p(V)}
             + \|g_{f_j}\|_{L^p(V)})
        \le 2C \|f_j\|_{\Np(V)}
\end{align*}
and hence
\begin{equation*}
\|u_j\|_{\Np(V)}
        \le \|w_j\|_{\Np(V)} + \|f_j\|_{\Np(V)}
       \le (1+2C) \|f_j\|_{\Np(V)}
\end{equation*}
is a bounded sequence in $\Np(V)$.
Since $u_j\to v$ and $w_j\to w$ q.e.\ in $X$, as $j\to\infty$,
it follows from Corollary~6.3 in~\cite{BBbook} that
$v\in \Np(V)$, $w\in \Np(X)$ and
\begin{equation}
    \|g_{v}\|_{L^p(V)}
    \le     \liminf_{j\to\infty} \|g_{u_j}\|_{L^p(V)}.
\label{eq-liminf-gj}
\end{equation}
From
$w=0$ q.e.\ in $X\setm V$ we see that $v-f_0\in\Np_0(V)$.
As $u_j\ge f_j\ge f_0$ q.e.\ in~$V$,
we have $v\ge f_0$ q.e.\ in $V$, and hence $v\in \K_{f_0,f_0}(V)$.

Next set $\phi_j=f_j + u_0-f_0
\in \K_{f_j,f_j}(V)$.
It follows that
\begin{equation}  \label{eq-uj-le-phij}
\int_V g_{u_j}^p\, d\mu\le \int_V g_{\phi_j}^p \, d\mu.
\end{equation}
Moreover, $\phi_j - u_0=f_j-f_0 \to 0$ in $\Np(V)$ and
in particular $g_{\phi_j} \to g_{u_0}$ in $\Lp(V)$, as $j \to \infty$.
Thus, using \eqref{eq-liminf-gj} and \eqref{eq-uj-le-phij},
\[
  \int_V g_{v}^p \, d\mu
      \le \liminf _{j\to \infty }\int_V g_{u_j}^p  \,d\mu
      \le \liminf _{j\to \infty }\int_V g_{\phi_j}^p \,d\mu
      = \int_V g_{u_0}^p \,d\mu.
\]
Therefore,
$v$ is also a solution of the $\K_{f_0,f_0}(V)$-obstacle problem,
and hence $v=u_0$ q.e., by the uniqueness part of
Theorem~\ref{thm-obstacle}.
\end{proof}

\section{Perron solutions}
\label{sect-fine-perron}

We are now ready to define and study Perron solutions on finely open sets.
We give four different definitions.
As we shall see they give almost the same solutions,
but this is far from immediate.
Recall that $\bdyp V$ denotes the fine boundary of $V$,
and that the fine regularizations $u_*$ and $u^*$
were defined in Section~\ref{sect-fine-cont}.

\begin{deff}  \label{def-P}
Given a function $f :\bdy_p V \to \eR$, let $\UU_f$
be the set of all
fine superminimizers $u\in\Np(V)$ such that
\begin{equation*} 
\fineliminf_{V \ni y \to x} u(y) \ge f(x) \quad \text{for all }
	x \in \bdy_p V.
\end{equation*}
Define the \emph{upper   $\sP$- and $\sR$-Perron solutions} of $f$ by
\begin{equation*} 
\uP f (x) = \inf_{u \in \UU_f}  u_*(x)
\quad \text{and} \quad \uR f (x) = \inf_{u \in \UU_f}  u^*(x), \qquad x \in V.
\end{equation*}

Similarly,  let $\UUt_f$ be the set of all
fine superminimizers $u\in\Np(V)$ such that
\begin{equation} \label{eq-UUt}
\fineliminf_{V \ni y \to x} u(y) \ge f(x) \quad \text{for q.e.\ }
	x \in \bdy_p V,
\end{equation}
and define the \emph{upper $\sS$- and $\sQ$-Perron solutions} as
\[ 
\uS f (x) = \inf_{u \in \UUt_f}  u_*(x)
\quad \text{and} \quad \uQ f (x) = \inf_{u \in \UUt_f}  u^*(x), \qquad x \in V.
\] 

The \emph{lower Perron solutions}
are defined analogously using suprema of fine subminimizers
in $\Np(V)$, or by letting
\[
\lP f:= - \uP (-f), \quad \lR f:= - \uR (-f), \quad \lS f:= - \uS (-f)
\quad \text{and} \quad \lQ f:= - \uQ(-f).
\]
\end{deff}

As usual, $\inf \emptyset = \infty$.
We shall call the fine regularizations
$u_*$ and $u^*$, appearing in the above infima,
\emph{admissible} in the definitions of the corresponding
Perron solutions.

Note that since functions in $\Np(V)$ are finite q.e.,
the upper Perron solutions are either
$<\infty$ q.e.\ or $\equiv\infty$.
Clearly,
\begin{equation}
\uS f \le \uP f \le \uR f
\quad \text{and} \quad
\uS f \le \uQ f \le \uR f \qquad \text{everywhere in $V$}.    \label{eq-diamond}
\end{equation}
We shall later see that under a very mild assumption,
all these solutions are equal q.e.\
and finely continuous q.e.
We postpone to Section~\ref{sect-Perron-fine-min} the proof that
Perron solutions are indeed fine minimizers,
under some growth conditions on the boundary data.

See Remark~\ref{rmk-not-lower-bdd} for why we do not need to
assume that functions in $\UU_f$ and $\UUt_f$ are bounded from below.
As
we only have
the comparison principle (Theorem~\ref{thm-qe-comparison-refined} below)
for functions in $\Np(V)$,
we instead require that functions in the upper classes $\UU_f$ and $\UUt_f$ belong
to $\Np(V)$.
This will be crucial when obtaining
Theorem~\ref{thm-lQ<=uQ} comparing
upper and lower 
$\sQ$-Perron solutions.

Since each function $u^*\in \UUt_f$ 
is finely upper semicontinuous, so are  $\uQ f$ and $\uR f$. Hence
\begin{equation} \label{eq-uuq}
  (\uQ f)_*\le (\uQ f)^*\le \uQ f
  \quad \text{and}\quad
  (\uR f)_*\le (\uR f)^*\le \uR f
  \qquad \text{everywhere in }V.
\end{equation}
Similarly, $\lQ f$ and $\lR f$ are finely lower semicontinuous and
\begin{equation} \label{eq-llq}
  (\lQ f)^*\ge (\lQ f)_*\ge \lQ f
  \quad \text{and}\quad
  (\lR f)^*\ge (\lR f)_*\ge \lR f
  \qquad \text{everywhere in }V.
\end{equation}
In particular, if $\uQ f= \lQ f$, then it is finely continuous.
We shall see in Proposition~\ref{prop-Q-reg}
and Lemma~\ref{lem-Qf=Rf-q.e.}
that
equality in~\eqref{eq-uuq} and~\eqref{eq-llq}
holds q.e.\ under a very mild assumption.

The use of finely \emph{upper semicontinuous} functions
in the definitions of $\uQ f$ and $\uR f$ is a key 
tool for obtaining
results for the Perron solutions,
see Remark~\ref{rmk-fine-min}.
Another practical reason for introducing $\uQ f$ is that
the inequality in \eqref{eq-UUt} is only required q.e., which
proves to be a useful tool.
This  also immediately implies that $\sQ$- and $S$-Perron solutions are
invariant under perturbations on sets of zero capacity,
i.e.\ that $\uQ f = \uQ f'$
and $\uS f = \uS f'$  whenever
$f=f'$ q.e.\ on $\bdy_p V$.

At the same time, the $P$-Perron solutions most closely
mimic the traditional definition of Perron solutions using
lsc-regularized \p-superharmonic functions
dominating the boundary data $f$ everywhere on the boundary.
This is the reason why we have reserved the letter $P$ for the
$P$-Perron solutions. The letter $Q$ in the $Q$-Perron solutions
comes from the q.e.-condition in their
definition. It is also closely related to the $Q$-Perron solutions
introduced in Bj\"orn--Bj\"orn--Shan\-mu\-ga\-lin\-gam~\cite{BBS2}.
The $S$-Perron solutions also share this property, but
for us the $Q$-Perron solutions are more important
because of Lemma~\ref{lem-uQ-defining-seq} and
Theorem~\ref{thm-lQ<=uQ} below.
Instead, $S$ comes from it being the smallest upper Perron solution.
Finally, $R$ was chosen since the upper and the lower $R$-solutions
are most relaxed (far away) with respect to each other,
  and as a tribute to Remak~\cite{remak}.
  In addition, all 
  four Perron solutions are closely related
to the Sobolev--Perron solutions
introduced in Bj\"orn--Bj\"orn--Sj\"odin~\cite{BBSjodin},
  see Remark~\ref{rmk-compare-Sjodin}.

It is the following countability lemma that makes
it possible to obtain some fundamental results for $\uQ f$.

\begin{lem} \label{lem-uQ-defining-seq}
Assume that $f:\bdyp V \to \eR$ is such that\/ $\UUt_f \ne \emptyset$.
Then there is a decreasing sequence $\{u_j\}_{j=1}^\infty$ of
functions in $\UUt_f$ such that
\[
  \uQ f(x)=\lim_{j \to \infty} u_j(x)
  \quad \text{for q.e.\ } x \in V.
\]
The same statement holds also for $\uR f$ and $u_j\in\UU_f$.
\end{lem}

\begin{proof}
Since each
function $u^*\in \UUt_f$, admissible for
$\uQ f$, is finely upper semicontinuous, so is $\uQ f$.
Moreover, for any $q\in \Q$,
\begin{equation} \label{eq-usc-union}
    \{x \in V : \uQ f(x) < q\} = \bigcup_{u \in \UUt_f} \{x \in V : u^*(x) < q\}
\end{equation}
is a union of finely open sets.
By the quasi-Lindel\"of principle (Theorem~\ref{thm-quasiLindelof}),
we can find a countable collection $\{u_{j,q}\}_{j=1}^\infty$ of functions in $\UUt_f$
and a set $E_q$ such that $\Cp(E_q)=0$ and
\[
    \{x \in V : \uQ f(x) < q\}
= E_q \cup \bigcup_{j=1}^\infty \{x \in V : u^*_{j,q}(x) < q\}.
\]
Letting $E=\bigcup_{q \in \Q} E_q$, we see that
\[
     \uQ f (x)= \inf_{\substack{q\in \Q \\ j=1,2,\ldots}} u^*_{j,q}(x)
     \quad \text{for all } x \in V \setm E.
\]
Next, we reorder the countable collection
\[
   \{u_{j,q}: q \in \Q  \text{ and } j=1,2,\ldots\}
\]
into a sequence $\{v_j\}_{j=1}^\infty$.
Then $u_j:=\min\{v_1,\ldots,v_j\} \in \UUt_f$
for each $j$ and the
sequence $\{u_j\}_{j=1}^\infty$ is decreasing.
Since $u^*_j\le \min\{v_1^*,\ldots,v_j^*\}\to \uQ f$ in
$V\setm E$, it follows that $u_j\to\uQ f$ q.e.\ in $V$.
The proof for $\uR f$ is the same.
\end{proof}

We can now show that the $\sQ$-Perron solutions
are finely continuous q.e.\ under a rather mild assumption.

\begin{prop} \label{prop-Q-reg}
Let $f: \bdyp V \to \eR$ and assume that $(\uQ f)_* > -\infty$ q.e.
Then
\begin{equation}  \label{eq-Q-reg}
(\uQ f)_* = (\uQ f)^* =\uQ f \quad \text{q.e.},
\end{equation}
and thus $\uQ f$ is finely continuous q.e.\ in $V$.
In particular, this holds if $\lQ f \not\equiv -\infty$.
\end{prop}

It follows from Theorem~1.4\,(b) in
Bj\"orn--Bj\"orn--Latvala~\cite{BBLat2} that
$\uQ f$ is also quasicontinuous on $V$, under the assumptions above.

\begin{proof}
We may assume that $\UUt_f\ne\emptyset$, since otherwise
$\uQ f \equiv \infty$ and there is nothing to prove.
Let $u_j\in \UUt_f$ be the decreasing sequence of functions provided by
Lemma~\ref{lem-uQ-defining-seq}.
Since $u_1\in\Np(V)$ is finite q.e.\ in~$V$, it
follows that for q.e.\ $x\in V$, there is
an integer $m$ such that
\[
x\in V_m:= \{y: (\uQ f)_*(y)+m > u^*_1(y)\}.
\]
Note that each $V_m$ is finely open
(because $(\uQ f)_*-u^*_1$
is finely lower semicontinuous) and that, by \eqref{eq-uuq},
\[
u_1-m < (\uQ f)_* \le \uQ f \le u_j \le u_1 \quad  \text{q.e.\ in }  V_m.
\]
Since each $u_j$ is a fine superminimizer in $V_m$, Theorem~\ref{thm-seq-fine-min-mod}
implies that so is $\uQ f$.
In particular, $\uQ f\in\Npploc(V_m)$ and it is thus finely continuous
q.e.\ in $V_m$.
This proves \eqref{eq-Q-reg} in $V_m$
and, by letting $m\to\infty$, also  in $V$.

Finally, if $v\in\Np(V)$ is admissible in the definition of 
 $\lQ f \not\equiv -\infty$, then
  \eqref{eq-llq} and the comparison in Theorem~\ref{thm-lQ<=uQ} below
  imply that
\[
(\uQ f)_* \ge (\lQ f)_* \ge \lQ f \ge v > -\infty \quad \text{q.e.}
\]
(Note that Proposition~\ref{prop-Q-reg} is not used
when proving Theorem~\ref{thm-lQ<=uQ} below.)
\end{proof}

An intriguing question is if $(\uQ f)^*=\uQ f$ everywhere.
In Proposition~\ref{prop-h-Qf-new} and
Theorem~\ref{thm-cont-Perron}  below
we show this in some special cases.
For Perron solutions on open sets, this is true
for arbitrary $f$, by Proposition~\ref{prop-open-set} below.
On the other hand, in the theory of balayage, which has many similarities with
the Perron method, it is known that one needs to regularize even on open sets.
We have found it helpful to use regularizations.

\section{Comparing Perron solutions}
\label{sect-comparison}

Our primary aim in this section is to deduce the following
fundamental inequalities between upper and lower  Perron solutions.

\begin{thm} \label{thm-lQ<=uQ}
Let $f: \bdyp V \to \eR$.
Then
\begin{equation*} 
\lQ f \le \lS f \le \uQ f 
\quad \text{everywhere in $V$}.
\end{equation*}
\end{thm}

A direct consequence of   \eqref{eq-diamond},
  Theorem~\ref{thm-lQ<=uQ} and the duality between upper and lower Perron solutions
  is that
\begin{equation} \label{eq-6.1alt}
\lR f \le \lQ f  \le \uS f \le \uQ f \le \uR f, \quad
\lR f  \le \lP f \le \uQ f
\quad \text{and} \quad 
\lQ f  \le \uP f \le \uR f
\end{equation}
everywhere in $V$.

We have not been able to prove
  that $\lP f \le \uP f$ and $\lS f \le \uS f$,
not even for bounded $f$.
A reason for this is that the inequality in our comparison principle
below holds only quasieverywhere.
This problem appears already for fine \p-supersolutions
in the Euclidean case, see Theorem~5.17 in Latvala~\cite{LatPhD}.
However when $p=n$ on unweighted $\R^n$,
the inequality in the comparison principle
was deduced everywhere in Corollary~7.13 in~\cite{LatPhD}.

The key to deducing Theorem~\ref{thm-lQ<=uQ}
is the following comparison principle.
The assumption $\Cp(X \setm V)>0$ from the beginning of
Section~\ref{sect-fine-min} is essential,
since otherwise \eqref{eq-a-liminf-ge-0-new} and~\eqref{eq-a-liminf-ge-0} are
void and the theorem fails,
see also Proposition~\ref{prop-cp-bdyp-U}.

\begin{thm} \label{thm-qe-comparison-refined}
\textup{(Comparison principle)}
Let $u \in \Np(V)$ be a fine superminimizer and
$v \in \Np(V)$ be a fine subminimizer.
Assume that for q.e.\ $x \in \bdyp V$, either
\begin{equation}   \label{eq-a-liminf-ge-0-new}
  \fineliminf_{V \ni y \to x}  u(y) \ge \finelimsup_{V \ni y \to x} v(y)
\end{equation}
or
\begin{equation}   \label{eq-a-liminf-ge-0}
  \fineliminf_{V \ni y \to x}  (u(y)-v(y)) \ge 0.
\end{equation}
Then $u \ge v$ q.e.\ in $V$, while
$u^* \ge v^* \ge v_*$ and $u^* \ge u_* \ge v_*$
everywhere in $V$.
\end{thm}

Note that \eqref{eq-a-liminf-ge-0-new} does not
imply \eqref{eq-a-liminf-ge-0} when the limits in
\eqref{eq-a-liminf-ge-0-new} are infinite.
The following lemma will be crucial for proving
Theorem~\ref{thm-qe-comparison-refined}.
It will also be used in the proof of Lemma~\ref{lem-Qf=Rf-q.e.}
and in Remark~\ref{rmk-not-lower-bdd}.

\begin{lem} \label{lem-UUt-UU-seq}
Let $u \in \Np(V)$ be a fine superminimizer and let $E\subset\bdy_p V$
with $\Cp(E)=0$.
Then there is a decreasing sequence $\{u_j\}_{j=1}^\infty$ of
finely lsc-regularized fine superminimizers in $\Np(V)$
which are bounded from below and
such that
\[
\lim_{j \to \infty} u_j=u   \text{ q.e.\ in }V
\quad \text{and} \quad
\lim_{V\ni y\to x} u_j(y) =\infty  \text{ for } x\in E \text{ and }
j=1,2,\ldots.
\]
If moreover $u \in \UUt_f$ for some $f: \bdyp V \to \eR$,
then in addition we can choose $u_j\in \UU_f$.
\end{lem}

\begin{proof}
By Corollary~1.3 in Bj\"orn--Bj\"orn--Shan\-mu\-ga\-lin\-gam~\cite{BBS5}
(or \cite[Theorem~5.31]{BBbook}),
$\Cp$ is an outer capacity.
Thus we can find a decreasing sequence
of open sets
$\{G_k\}_{k=1}^\infty$
such that $E \subset G_k$ and $\Cp(G_k) < 2^{-kp}$, $k=1,2,\ldots$\,.
By Lemma~10.17 in~\cite{BBbook} there is a decreasing sequence of nonnegative
functions $\{\psi_j\}_{j=1}^\infty$
such that $\| \psi_j \|_{\Np(X)} < 2^{-j}$ and
$\psi_j \ge k-j$ in $G_k$ whenever $k > j$.

Let $v_j = \psi_j+\max\{u,-j\}$ and let $u_j$ be the finely lsc-regularized
solution of the $\K_{v_j,v_j}(V)$-obstacle problem.
Then $u_j$ is a fine superminimizer,
by Theorem~\ref{thm-obstacle}.
Moreover, $u_j \ge -j$ and it is thus bounded from below.
If $k >j$, then
\[
v_j \ge \psi_j-j \ge k-2j \quad \text{in } G_{k} \cap V,
\]
and hence $u_j \ge k-2j$ in $G_{k} \cap V$.
Thus for each $j$,
\begin{equation}    \label{eq-infty-on-E}
\liminf_{V\ni y\to x} u_j(y)
=\infty \quad \text{for } x\in E.
\end{equation}
Since $v_j \searrow u$ in $\Np(V)$, and
$u$ is a solution of the $\K_{u,u}(V)$-obstacle problem
(by Lemma~\ref{lem-super-obst}), it follows from
Proposition~\ref{prop-3.2} that
$\{u_{j}\}_{j=1}^\infty$ decreases q.e.\ to~$u$.

Finally, if $u \in \UUt_f$ then
we let $\Et=E \cup \{x: \fineliminf_{V\ni y\to x} u(y) < f(x)\}$
and construct $u_j$ as above with $\Et$ instead of $E$.
(Note that $\Cp(\Et)=0$, as $u \in \UUt_f$.)
Thus
\[
\fineliminf_{V\ni y\to x} u_j(y)
\ge \fineliminf_{V\ni y\to x} u(y)
\ge f(x)
\quad \text{for } x\in\bdy_p V\setm \Et,
\]
which together with~\eqref{eq-infty-on-E} (with $E$ replaced by $\Et$)
shows that $u_j \in \UU_f$.
\end{proof}

\begin{proof}[Proof of Theorem~\ref{thm-qe-comparison-refined}]
Assume to start with that $u$ is bounded from below and $v$ from above.
Then \eqref{eq-a-liminf-ge-0-new}$\imp$\eqref{eq-a-liminf-ge-0}
and
Proposition~\ref{prop-Np0} shows that $\max\{v-u,0\} \in \Np_0(V)$.

Now let  $\vt$ be a solution of the Dirichlet
$\K_{-\infty,v}(V)$-obstacle problem.
Note that $-\vt$ is a solution of the Dirichlet
$\K_{-\infty,-v}(V)$-obstacle problem.
By Lemma~\ref{lem-super-obst},
$u$ is a solution of the $\K_{u,u}(V)$-obstacle problem,
and $-v$ is a solution of the $\K_{-v,-v}(V)$-obstacle problem.
By the comparison principle (Lemma~\ref{lem-obst-le}),
we therefore see that
$-\vt \le -v$ and $\vt \le u$ q.e.\ in $V$ and thus
\[
v\le \vt \le u  \quad \text{q.e.\ in } V.
\]

Next, for general $u$ and $v$, Lemma~\ref{lem-UUt-UU-seq}
provides us with a decreasing sequence $\{u_j\}_{j=1}^\infty$
of fine superminimizers bounded from below
such that $\lim_{j \to \infty} u_j(x)=u(x)$ q.e.\ in $V$.
In particular, $u_j \ge u$ q.e.\ in $V$.
Similarly, there is an
increasing sequence $\{v_j\}_{j=1}^\infty$
of fine subminimizers bounded from above
such that $\lim_{j \to \infty} v_j(x)=v(x)$ q.e.\ in $V$ and $v_j \le v$
q.e.\  in $V$.

Applying the already proved bounded case to $u_j$ and $v_j$,
we get $u_j \ge v_j$ q.e.\ in~$V$.
Letting $j \to \infty$ shows that $u \ge v$ q.e.\ in $V$
also in this case.
The statement about fine regularizations now follows immediately.
\end{proof}

We are now ready to prove Theorem~\ref{thm-lQ<=uQ}.

\begin{proof}[Proof of Theorem~\ref{thm-lQ<=uQ}]
Let $u^*$ and $v^*$ be admissible in the
definitions of $\uQ f$ and $\lS f$, respectively.
Then
\[
     \fineliminf_{V \ni y \to x} u^*(y) \ge f(x) \ge \finelimsup_{V \ni y \to x} v^*(y)
  \quad \text{for q.e.\ } x \in \bdyp V.
\]
It thus follows from
the comparison principle (Theorem~\ref{thm-qe-comparison-refined})
that $u^* \ge v^*$ everywhere in $V$.
Taking infimum resp.\ supremum over all admissible $u^*$ and $v^*$
yields that $\uQ f \ge \lS f$ in $V$.
By duality, the other inequality follows from~\eqref{eq-diamond}.
\end{proof}

We conclude this section by comparing the four definitions of Perron solutions
and prove that they coincide q.e.

\begin{thm}  \label{thm-all-equal}
Let $f: \bdyp V \to \eR$. Then
\begin{equation}   \label{eq-ineq-Q-S}
\lQ f  \le (\lQ f)_* \le (\uQ f)_* \le \uS f
\quad \text{everywhere in $V$}
\end{equation}
and
\begin{equation}   \label{eq-eq-SPQR}
(\uS f)_* = (\uP f)_* = (\uQ f)_* = (\uR f)_*
\quad \text{everywhere in $V$}.
\end{equation}
If moreover $(\uQ f)_*>-\infty$ q.e., then also
\[
\uS f = \uP f = \uQ f = \uR f \quad \text{q.e.\ in $V$} 
\]
and these four upper Perron solutions are finely continuous q.e.
\end{thm}

Note that by \eqref{eq-eq-SPQR} the condition $(\uQ f)_* > -\infty$ q.e.\
in the last part of Theorem~\ref{thm-all-equal},
as well as in Proposition~\ref{prop-Q-reg}, can equivalently
be expressed using any of the lsc-regularized upper Perron solutions.
In order to prove Theorem~\ref{thm-all-equal}
we will need the following two lemmas.

\begin{lem} \label{lem-S-Q-new}
Let $f: \bdyp V \to \eR$.
Then
\begin{equation*}
(\uQ f)_* \le \uS f \quad \text{and}  \quad (\uR f)_* \le \uP f
\qquad \text{everywhere in $V$}.
\end{equation*}
\end{lem}

\begin{proof}
Let $u\in \UUt_f$.
Then $u^*$ is admissible in the definition of $\uQ f$
and hence
$(\uQ f)_* \le (u^*)_* =u_*$ in $V$,
by Remark~\ref{rmk-reg} together with the fact that $u$
is finely continuous q.e.\ in $V$.
Since $u_*$ is admissible in the definition of $\uS f$,
taking infimum over all $u\in\UUt_f$ shows that $(\uQ f)_* \le \uS f$.

The inequality $(\uR f)_* \le \uP f$ is
shown similarly
by considering $u \in \UU_f$.
\end{proof}

\begin{lem}   \label{lem-Qf=Rf-q.e.}
Let $f: \bdyp V \to \eR$.
Then $\uQ f= \uR f$ q.e.\ in $V$.
\end{lem}

\begin{proof}
We may assume that $\UUt_f\ne\emptyset$.
Lemma~\ref{lem-uQ-defining-seq} provides us with a decreasing
sequence $\{u_j\}_{j=1}^\infty$ of functions in $\UUt_f$ such that
\[
\uQ f(x)=\lim_{j \to \infty} u_j(x)
  \quad \text{for q.e.\ } x \in V.
\]
For each $j=1,2,\ldots$\,, let $u_{j,k}\in \UU_f$ be the decreasing
sequence of functions
provided by Lemma~\ref{lem-UUt-UU-seq} so that
\[
u_j(x) = \lim_{k \to \infty} u_{j,k}(x)
\quad \text{for q.e.\ $x\in V$.}
\]
Then
\[
v_k := \min_{j\le k} u_{j,k}  \to \uQ f
\quad \text{q.e.\ in }V, \text{ as } k \to \infty.
\]
Since  $v_k^*$ are admissible in the definition of $\uR f$
and $v_k^*=v_k$ q.e., letting
$k\to\infty$ shows that $\uQ f\ge \uR f$ q.e.
The reverse inequality $\uQ f\le \uR f$ holds everywhere and is immediate
from the definition.
\end{proof}

\begin{proof}[Proof of Theorem~\ref{thm-all-equal}]
The inequalities~\eqref{eq-ineq-Q-S} follow directly from
\eqref{eq-llq},  Theorem~\ref{thm-lQ<=uQ} and Lemma~\ref{lem-S-Q-new}.
Next, Lemmas~\ref{lem-S-Q-new} and~\ref{lem-Qf=Rf-q.e.}
together with \eqref{eq-diamond}
yield
\begin{equation}   \label{eq-before-lsc-reg}
(\uQ f)_* \le \uS f \le \uP f
\le \uR f = \uQ f
\quad \text{q.e.\  in $V$}.
\end{equation}
Applying the fine lsc-regularization to this chain of inequalities,
together with Lemma~\ref{lem-u*-fine-reg}, gives
the equalities~\eqref{eq-eq-SPQR}.
The last statement 
follows from
\eqref{eq-before-lsc-reg} and Proposition~\ref{prop-Q-reg}.
\end{proof}

\begin{remark} \label{rmk-not-lower-bdd}
Traditionally, Perron solutions (on open sets) are defined using
functions in the upper class that are bounded from below.
This ``lower boundedness'' assumption was introduced by Brelot~\cite{brelot},
who was the first to study Perron solutions for unbounded functions;
see \cite[p.~146]{brelot} both for the definition and a remark
explaining why it is essential.
The first papers on Perron solutions,
by Perron~\cite{Perron23} and Remak~\cite{remak}, only dealt
with continuous (and thus bounded) boundary data $f$,
for which this assumption is redundant.

As we have seen, we do not need this assumption.
The reason for this is that
our upper and lower classes only contain Sobolev functions, for which
the comparison
principle (Theorem~\ref{thm-qe-comparison-refined}) holds
without any boundedness assumption.
On the other hand, even with this assumption,
essentially all of the theory
developed here would be true,
including  the results in this section.
Let us denote such
upper Perron solutions by
$\uShat f$,   $\uPhat f$, $\uQhat f$ and $\uRhat f$.
By Lemmas~\ref{lem-uQ-defining-seq} and~\ref{lem-UUt-UU-seq}
it follows (essentially as in the proof of Lemma~\ref{lem-Qf=Rf-q.e.})
that 
$\uRhat f \le \uQ f$ q.e.
Using also the trivial inequalities
$\uQ f \le \uQhat f \le \uRhat f$ gives
\begin{equation} \label{eq-hat}
  \uQhat f = \uRhat f = \uQ f = \uR f
  \quad \text{q.e.}
\end{equation}  
The $\widehat{\hphantom{P}}$-version of Theorem~\ref{thm-all-equal} then yields
that 
\begin{equation}    \label{eq-hats=bar}
\uShat f = \uPhat f = \uQhat f= \uRhat f= \uQ f \quad \text{q.e.},
\end {equation}
provided that
$(\uQhat f)_*>-\infty$ q.e.\
(or equivalently
$(\uQ f)_*>-\infty$ q.e.).
See also Remark~\ref{rmk-not-lower-bdd-open} below.

The only real difference after adding
the ``lower boundedness'' assumption
would be in Proposition~\ref{prop-h-Qf-new} below,
where $h^*$ with $\inf h^* =-\infty$
is not
admissible in the definition
of $\uQhat f$
and thus we would not be able to conclude that $\uQhat f =h^*$
and $\uShat f =h_*$ everywhere for such $h$.
It would still follow that $\uQhat f= \uQ f= h^*$ 
$\uShat f= \uS f= h_*$ q.e., by \eqref{eq-hat} and \eqref{eq-hats=bar}.
\end{remark}

\section{Perron solutions are fine minimizers}
\label{sect-Perron-fine-min}

Our next aim is to show that the Perron solutions
are indeed fine minimizers, under rather general assumptions
on the boundary data.

\begin{thm} \label{thm-uQ-finemin}
Let $f: \bdyp V \to \eR$.
Assume that one of the following conditions
holds\/\textup{:}
\begin{enumerate}
\item \label{g-ba}
$\uQ f\in \Npploc(V)$\textup{;}
\item \label{g-bb}
there is $M>0$ and $u \in \Np(V)$ such that
$|\uQ f -u| \le M$ a.e.\ in $V$\textup{;}
\item \label{g-bc}
there is $M>0$ and $u \in \Np(V)$
such that
\begin{equation*}
\finelimsup_{V\ni y\to x}u(y)- M \le f(x) \le \fineliminf_{V\ni y\to x}u(y)+ M
  \quad \text{for q.e.\ } x\in\bdy_p V.
\end{equation*}
\end{enumerate}
Then $\uQ f$ is a fine minimizer. In particular, this holds if $f$ is q.e.\ bounded.
\end{thm}

As a direct consequence we obtain the following characterization.

\begin{cor}
Let $f: \bdyp V \to \eR$.
Then $\uQ f$ is a fine minimizer if and only if $\uQ f\in \Npploc(V)$.
\end{cor}

A weaker necessary condition for $\uQ f$ to be a fine minimizer is that
$\uQ f$ is finite q.e. We do not know if this condition is sufficient.

To prove Theorem~\ref{thm-uQ-finemin},
we will need the following lemma.

\begin{lem} \label{lem-Poisson-all-V}
For $f: \bdyp V \to \eR$ and $u \in \UUt_f$, 
let $h$ be a
solution of the Dirichlet $\K_{-\infty,u}(V)$-obstacle problem.
Then $h \in \UUt_f$, $h \le u$ q.e.\
in $V$ and $h$ is a fine minimizer.
\end{lem}

\begin{proof}
By Theorem~\ref{thm-obstacle}, $h$ is a fine minimizer.
Since $u$ is a solution of the $\K_{u,u}(V)$-obstacle
problem, by Lemma~\ref{lem-super-obst},
the comparison principle (Lemma~\ref{lem-obst-le})
implies that  $h \le u$ q.e.\ in $V$.
By definition, $h-u \in\Np_0(V)$,
and thus by Proposition~\ref{prop-Np0}, also $h \in \UUt_f$.
\end{proof}

\begin{proof}[Proof of Theorem~\ref{thm-uQ-finemin}]
By Lemma~\ref{lem-uQ-defining-seq},
there is a decreasing sequence $\{u_j\}_{j=1}^\infty$ of
functions in $\UUt_f$ such that
\[
  \uQ f(x)=\lim_{j \to \infty} u_j(x)
  \quad \text{for q.e.\ } x \in V.
\]
If \ref{g-ba} or \ref{g-bb} in the statement of the theorem
holds, then replace
each $u_j$ by a fine minimizer $h_j\in \UUt_f$ provided by
Lemma~\ref{lem-Poisson-all-V}, so that $h_j\to\uQ f$
q.e.\ in $V$ as $j \to \infty$.
It follows from the comparison principle (Lemma~\ref{lem-obst-le}) that
(after redefinitions on sets of capacity zero)
$\{h_j\}_{j=1}^{\infty}$ is also a decreasing sequence.
Theorem~\ref{thm-seq-fine-min-mod} then shows that $\uQ f$
is a fine minimizer.

Assume next that \ref{g-bc} holds and let $h\in \Np(V)$ be
a solution of the Dirichlet $\K_{-\infty,u}(V)$-obstacle problem.
Since $h-u\in\Np_0(V)$,
Proposition~\ref{prop-Np0} shows that
$\finelim_{V\ni y\to x}(h-u)(y)=0$ for q.e.\ $x\in\bdy_p V$,
and hence
\begin{equation}   \label{eq-f-le-h'+M}
f(x) \le \fineliminf_{V\ni y\to x}u(y)+ M = \fineliminf_{V\ni y\to x}h(y)+ M
\quad \text{for q.e.\ $x\in\bdy_p V$.}
\end{equation}
Thus $v_j:=\min\{u_j,h+M\}\in\UUt_f$.

Let $h_j\in\UUt_f$ be a fine minimizer provided for $v_j$
by Lemma~\ref{lem-Poisson-all-V}.
As in~\eqref{eq-f-le-h'+M}, we have using~\ref{g-bc} that
\[
\finelimsup_{V\ni y\to x} h(y)-M \le f(x) \le \fineliminf_{V\ni y\to x} h_j(y)
\quad \text{for q.e.\ $x\in\bdy_p V$.}
\]
The comparison principle (Theorem~\ref{thm-qe-comparison-refined}) then implies that
\[
h-M \le h_j \le v_j \le h+M \quad \text{q.e.\ in $V$.}
\]
Since $h_j\to\uQ f$ q.e.\ in $V$ as $j \to \infty$,
Theorem~\ref{thm-seq-fine-min-mod} concludes the proof.
\end{proof}

\begin{remark} \label{rmk-h}
Another consequence of Lemma~\ref{lem-Poisson-all-V}
and Proposition~\ref{prop-Np0} is that 
\[
  \uQ f = \inf\{h^* \in \UUt_f : h \text{ is a fine minimizer}\}
  = \inf_u H^*u
\quad \text{everywhere in $V$,}
  \]
where the second infimum is taken over all $u \in \Np(V)$ such that
\eqref{eq-UUt} holds, and $H^*u$ is the usc-regularized  solution of
the Dirichlet $\K_{-\infty,u}(V)$-obstacle problem.
Similar identities hold for $\uS f$.
On the other hand, it is far from clear whether $\uP f$ and $\uR f$
can be expressed in such terms.
\end{remark}  

\begin{remark} \label{rmk-fine-min}
Let us compare the proof of
Theorem~\ref{thm-uQ-finemin} with the more traditional proofs for
Perron solutions on open sets,
as in
Heinonen--Kilpel\"ainen--Martio~\cite[Theorem~9.2]{HeKiMa} and
Bj\"orn--Bj\"orn--Shan\-mu\-ga\-lin\-gam~\cite[Theorem~4.1]{BBS2}
(or \cite[Theorem~10.10]{BBbook}).
A major idea in those proofs is to use convergence on
a countable dense subset of the open set $V$.
(In \cite{HeKiMa} this appears in the proof of
Choquet's topological lemma~\cite[Lemma~8.3]{HeKiMa}.)
In many cases, countable subsets of $X$ have zero capacity and are
therefore not seen by the fine topology.
Hence there is no possibility to obtain general results using a countable
finely dense subset of a finely open $V$.
Instead,
we use the countable dense subset $\Q$ on the target
side via Lemma~\ref{lem-uQ-defining-seq} and obtain convergence q.e.\
in $V$.
For this, the use of finely usc-regularized functions is essential as it
makes it possible to deduce the crucial identity~\eqref{eq-usc-union}
in Lemma~\ref{lem-uQ-defining-seq} with finely open sets
in the right-hand side.

Another ingredient in the 
traditional proofs is the
Poisson modification on an exhaustion by compactly
contained open subsets of $V$, together with Harnack's
convergence principle.
For fine minimizers, we have only weaker convergence theorems
which require the additional assumptions in Theorem~\ref{thm-uQ-finemin}.
At the same time, since we define the $Q$-Perron solutions
using functions from
$\Np(V)$ and 
require the boundary
inequality in~\eqref{eq-UUt} only for q.e.\ $x\in\bdy_p V$,
the Poisson modification in the proof of Theorem~\ref{thm-uQ-finemin}
can be taken with respect to all of $V$ through Lemma~\ref{lem-Poisson-all-V}.
\end{remark}

\section{Resolutivity of Sobolev and continuous boundary values}
\label{sect-continuity}

Next we study \emph{resolutivity},
i.e.\ when the upper and lower Perron solutions agree,
at least q.e.
As we shall see, the $\sQ$-Perron solutions seem to behave
slightly better than the other solutions.
In combination with 
Theorem~\ref{thm-all-equal}, they provide resolutivity and
invariance results for the other Perron solutions as well.

\begin{prop}  \label{prop-h-Qf-new}
Let $f : \clVp \to \eR$.
Assume that $f\in\Np(V)$ and
\begin{equation} \label{eq-prop-h-Qf}
f(x)=\finelim_{V\ni y\to x}f(y)
\quad \text{for q.e.\ $x\in\bdy_p V$}.
\end{equation}
In particular, this holds if $f\in\Np(X)$.

Let $h$ be a solution of the
Dirichlet $\K_{-\infty,f}(V)$-obstacle problem.
Then
\begin{equation}
\uQ f = h^* =\lS f \quad \text{and} \quad \lQ f = h_* =\uS f
\qquad \text{everywhere in } V
\label{eq-Qf=h=Sf}
\end{equation}
and $\lQ f = \uQ f$ q.e.\ in $V$.
Moreover,
\begin{equation}
(\lQ f)^* = (\uQ f)^* = \uQ f  \quad \text{and} \quad
(\uQ f)_* = (\lQ f)_* = \lQ f
\qquad \text{everywhere in } V.
\label{eq-Qf=Qf*}
\end{equation}
\end{prop}

\begin{proof}
Proposition~\ref{prop-Np0}, Theorem~\ref{thm-obstacle}, and the assumptions on $f$
imply that $h^*$ is admissible in the definition of
both $\uQ f$ and $\lS f$.
Hence by Theorem~\ref{thm-lQ<=uQ},
\[
\uQ f \le h^* \le \lS f \le \uQ f \quad \text{everywhere in } V,
\]
which proves the first statement in~\eqref{eq-Qf=h=Sf}.
The second statement in~\eqref{eq-Qf=h=Sf} is shown similarly,
or by applying the first one to $-f$.
The equality $\lQ f = \uQ f$ q.e.\ then follows immediately
from the fact that $h^*=h_*$ q.e.
Applying the  fine regularizations, together with
Lemma~\ref{lem-u*-fine-reg}, now gives
\[
(\lQ f)^* = (\uQ f)^* = h^{**} = h^* = \uQ f
\]
as well as the other identities in~\eqref{eq-Qf=Qf*}.

Finally, if $f \in \Np(X)$, then $f$ is finely continuous q.e.\
in $X$ and thus \eqref{eq-prop-h-Qf} holds.
\end{proof}

The following example shows that
the fine limit in~\eqref{eq-prop-h-Qf}
need not exist for a large part of $\bdy_p V$ when $f \in \Np(V)$.

\begin{example}
Let (using complex notation)
\[
V= \{z=re^{i\theta}: 0<r<1 \text{ and } 0<\theta<2\pi \}
\]
be the slit disc in $\R^2$.
Then $f(re^{i\theta}):= r\theta \in \Np(V)$ for all $p>1$,
but the fine limits
\[
\finelim_{V\ni y\to z} f(z)
\]
do not exist for any $z$ in the slit, but for the tip $0$.
This also shows that condition~\ref{g-bc} in Theorem~\ref{thm-uQ-finemin}
would be less general if
$\finelimsup$ and $\fineliminf$ were replaced by $\finelim$.
\end{example}  

We shall now see that under additional assumptions on
the boundary data, the $Q$-Perron solutions are finely continuous.

\begin{prop}  \label{prop-h-fine-cont}
Let
\[
  V_0 = \{z \in V : \Cp(B(z,r) \cap \bdyp V)>0 \text{ for all } r >0\}.
\]
Also let $f : \clVp \to \eR$ be such that $f\in\Np(V)$,
\eqref{eq-prop-h-Qf} holds and
\begin{equation}   \label{eq-cplim-ex}
\cplimalt_{\bdy_p V\ni x\to z} f(x) \quad \text{exists for all
  $z\in V_0$.}
\end{equation}
In particular, this holds if $f=f_0$ q.e.\
in $V$ for some $f_0 \in \Lip(\clVp)$.

Then $h_* = h^*$ is finely continuous and
\begin{equation}
\lS f = \lQ f = h_* = h^* = \uQ f = \uS f
   \quad \text{everywhere in } V.
  \label{eq-k1-new}
\end{equation}
\end{prop}

Here $\cplim$ is taken with respect
to the metric topology from $X$ and up to sets of zero capacity,
as in~\cite[Section~7]{BBLat4}.
Note that the limit in~\eqref{eq-cplim-ex} does not make sense
for $z\in V\setm V_0$.
The simplest example showing that $V_0$ can be nonempty is
  perhaps letting $V=\Om \cup \{0\}$, where $\Om$ is
  the complement of the Lebesgue spine in $\R^3$ as in Example~13.4
  in~\cite{BBbook}.

\begin{proof}
The assumption \eqref{eq-prop-h-Qf} implies that 
the extensions 
\[
f_*(x) = \fineliminf_{V \ni y \to x} f(y)
\quad \text{and} \quad
f^*(x) = \finelimsup_{V \ni y \to x} f(y)
\qquad \text{for } x \in \bdyp V
\]
satisfy $f_*=f=f^*$ q.e.\ on $\bdyp V$.
Hence, as \eqref{eq-cplim-ex} holds, 
all $z\in V_0$ satisfy
condition~(b) in \cite[Theorem~7.5]{BBLat4}.
On the other hand, for $z\in V\setm V_0$,
condition~(c) in \cite[Theorem~7.5]{BBLat4} holds.
Theorem~7.5 in \cite{BBLat4} then implies that $h_*=h^*$
is finely continuous in $V$.
The remaining conclusions in~\eqref{eq-k1-new} then immediately follow
from~\eqref{eq-Qf=h=Sf}.
\end{proof}

Our next aim is to prove the following resolutivity
and invariance result.
Here $\Cunif(\bdyp V)$ is the space of
uniformly continuous functions on $\bdyp V$,
with respect to the metric topology.

\begin{thm}\label{thm-cont-Perron}
Let $f \in \Cunif(\bdyp V)$ and $k:\bdy_p V\to\eR$ be a function which vanishes q.e.
Then
\begin{equation*}
\lS f = \uS f = \lQ f = \uQ f = \lS (f+k) = \uS (f+k) = \lQ (f+k) = \uQ (f+k)
\end{equation*}
everywhere in $V$,
and this function  is finely \p-harmonic,
i.e.\ a finely continuous fine minimizer.
\end{thm}

A natural question is if some resolutivity can be shown
for finely continuous functions
on the boundary. However, even if $V\subset\R^n$ is open,
there are no such results available
in the nonlinear theory.
Another question is if it would be enough to require
$f \in C(\bdyp V)$.
Due to the possible noncompactness
of the fine boundary with
respect to the metric topology, this is not
equivalent to requiring that $f \in \Cunif(\bdyp V)$.
Resolutivity for functions in $C(\bdyp V)$
is unknown even for open $V$, cf.\
the discussion in Bj\"orn~\cite{ABcomb}.

\begin{proof}[Proof of Theorem~\ref{thm-cont-Perron}]
Since $X$ is proper and $V$ is bounded,  $\bdy_p V$ is totally bounded.
Using that $f$ is uniformly continuous, we can therefore for
each $j=1,2,\ldots$\,, find $\de_j>0$ and finitely many
balls $B_{i,j}:=B(x_{i,j},\de_j)$ covering $\bdy_p V$, such that
\[
\osc_{2B_{i,j}\cap \bdy_p V} f <1/j.
\]
Let 
\[
\psi_{i,j}(x) =
\max\{0,\min\{1,2-d(x,x_{i,j})/\de_j\}\}
\quad \text{and} \quad
\eta_{i,j} = \frac{\psi_{i,j}}{\max\{1,\sum_{i} \psi_{i,j}\}}
\]
be a Lipschitz partition of unity on $\bigcup_i B_{i,j}$
subordinate to $2B_{i,j}$.
Then 
\[
f_j:= \sum_i f(x_{i,j}) \eta_{i,j}\in \Lipc(X)
\]
and 
\[
|f(x)- f_j(x)| \le \sum_i |f(x)-f(x_{i,j})| \eta_{i,j}(x) < 1/j
\quad \text{for } x\in \bdyp V.
\]
Hence, by Proposition~\ref{prop-h-fine-cont}, applied to $f_j+1/j$,
\[
    \uQ f \le \uQ (f_j + 1/j) = \lQ (f_j +1/j)
    \le 2/j + \lQ f
    \quad \textrm{in }V.
\]
Letting $j \to \infty$ shows that $\uQ f \le \lQ f$ in $V$.
As also the converse inequality holds,
by Theorem~\ref{thm-lQ<=uQ}, we see that $\uQ f = \lQ f$.
The equalities $\lQ (f+k)=\lQ f$ and $\uQ f=\uQ (f+k)$
follow directly from the
definition and the fact that $k=0$ q.e.
Finally, Theorem~\ref{thm-lQ<=uQ} and \eqref{eq-6.1alt}
include also $\lS f$, $\uS f$,
$\lS(f+k)$ and $\uS(f+k)$ in the equalities.
\end{proof}

\section{Perron solutions on open sets}
\label{sect-open}

\emph{In addition to the standing assumptions from the beginning of
Section~\ref{sect-fine-cont},
we assume in this section that $V$ is a bounded open set with $\Cp(X \setm V)>0$.}

\medskip

In this final section we will show that our four upper Perron solutions
coincide if $V$ is open.
Note that for open sets $V$, the definitions of the  $P$-, $Q$- and $S$-Perron solutions
in \cite{BBbook}, \cite{BBS2} and Bj\"orn--Bj\"orn--Sj\"odin~\cite{BBSjodin}
are different from the ones considered here.
One  major difference is that our Perron solutions use
fine limits on the fine boundary $\bdyp V$
instead of
ordinary limits on the full metric boundary $\bdy V$,
  see Remark~\ref{rmk-compare-Sjodin} below.

\begin{prop} \label{prop-open-set}
Assume that $V$ is open and let $f : \clVp \to \eR$.
Then
\[
\uS f = \uP f = \uQ f = \uR f \quad \text{everywhere in $V$}.
\]
Moreover, this upper Perron solution
is in each component of $V$ either identically $\pm \infty$ or
\p-harmonic, i.e.\ a continuous minimizer.
\end{prop}

\begin{proof} We may assume that $\UUt_f \ne \emptyset$, as otherwise
    $\uS f=\uP f =\uQ f = \uR f= \infty$ everywhere in $V$. Let $u \in \UUt_f$.
As $V$ is open, $u$ is a standard superminimizer
(as e.g.\ in~\cite{BBbook}), by
Corollary~5.6 in~\cite{BBLat4}.
It thus follows from
Theorem~5.1 and Proposition~7.4 in Kinnunen--Martio~\cite{KiMa02}
(or Theorem~8.22 and Proposition~9.4 in~\cite{BBbook})
that $u$ has a \p-superharmonic representative, which is finely continuous
by   Bj\"orn~\cite[Theorem~4.4]{JB-pfine} or
Korte~\cite[Theorem~4.3]{korte08} (or \cite[Theorem~11.38]{BBbook}), i.e.\
$u^* = u_*$.
As this holds for all $u \in \UUt_f$, and in particular for all $u \in \UU_f$,
we get immediately from the definitions that
\[
     \uS f = \uQ f
\quad \text{and} \quad
     \uP f = \uR f
\qquad \text{everywhere in } V.
\]
Similarly, Lemma~\ref{lem-Poisson-all-V} implies that $\uQ f$ is
a pointwise infimum of continuous minimizers,
  and is thus
upper semicontinuous in $V$.

One can now proceed essentially verbatim as in
Bj\"orn--Bj\"orn--Shan\-mu\-ga\-lin\-gam~\cite[Theorem~4.1]{BBS2}
(or \cite[Theorem~10.10]{BBbook}) to show that the Perron
solutions are \p-harmonic.
Alternatively one can argue as follows:

For $\uR f$, Lemma~\ref{lem-uQ-defining-seq} provides us with
a decreasing sequence
of \p-superharmonic functions $u_j\in \UU_f$ such
that $u_j\to \uR f$ q.e.\ in $V$.
For each $j$, let $\ut_j\in \UU_f$ be the
Poisson modification of $u_j$ in an open set $\Om\Subset V$,
as in \cite[Theorem~9.44]{BBbook}.
Then
$\ut_j \le u_j$  and $\{\ut_j\}_{j=1}^\infty$ is also a decreasing sequence.
Let $\ut=\lim_{j\to\infty} \ut_j$ in $V$.
It follows that
$\ut\ge \uR f \ge \uQ f$ everywhere in $V$.
By Lemma~\ref{lem-Qf=Rf-q.e.}, $\ut = \uR f= \uQ f$  in $V\setm E$
  for some set $E$ with $\Cp(E)=0$.
  Note that $V\setm E$ is dense in $V$.

Since each $\ut_j$ is \p-harmonic in $\Om$, \cite[Corollary~9.38]{BBbook}
implies that in every component of $\Om$, either
$\ut$ is \p-harmonic or $\ut \equiv -\infty$.
We thus have in each such component, using the upper semicontinuity of
$\uQ f$ and the continuity of $\ut$, that
\[
\ut(x) \ge \uR f(x) \ge \uQ f(x) \ge \limsup_{\Om\setm E\ni z\to x} \uQ f(z)
= \limsup_{\Om\setm E\ni z\to x} \ut(z) = \ut(x),
  \quad x \in V.
\]
Thus $\uR f = \uQ f =\ut$ is \p-harmonic, or identically $-\infty$,
in each component of $\Om$ and as $\Om\Subset V$
was arbitrary, also in each component of $V$.
\end{proof}

\begin{remark} \label{rmk-not-lower-bdd-open}
Consider the upper Perron solutions 
  $\uShat f$,   $\uPhat f$, $\uQhat f$ and $\uRhat f$
  introduced in Remark~\ref{rmk-not-lower-bdd}.
The proof of Proposition~\ref{prop-open-set} applies
equally well to them, showing that
$\uShat f = \uPhat f = \uQhat f = \uRhat f$ everywhere in $V$,
and that this function is $\eR$-continuous.
As shown
in Remark~\ref{rmk-not-lower-bdd}, 
$\uQhat f=\uQ f$ q.e.,
but since
both functions are $\eR$-continuous they coincide everywhere.
Hence, for open $V$,
\[
\uShat f = \uPhat f = \uQhat f = \uRhat f
= \uS f = \uP f = \uQ f = \uR f \quad \text{everywhere in $V$}.
\]
\end{remark}

\begin{remark} \label{rmk-compare-Sjodin}
  Upper Perron solutions on bounded open sets in metric spaces,
  defined using 
  lower bounded \p-superharmonic functions satisfying 
\begin{equation}  \label{eq-liminf}
\liminf_{V\ni y\to x} u(y) \ge f(x)
\end{equation}
on the whole boundary $\bdy V$ rather than  $\bdy_p V$,
were first 
studied in
Bj\"orn--Bj\"orn--Shan\-mu\-ga\-lin\-gam~\cite{BBS2}.
Such solutions are traditionally denoted $\uP f$, but to avoid confusion with our 
definitions, we use $\tP f$ for them.
Similarly, in Bj\"orn--Bj\"orn--Sj\"odin~\cite{BBSjodin},
upper Sobolev--Perron solutions (denoted $\uS f$ therein)
were defined by requiring in addition that 
the admissible \p-superharmonic functions $u$ have finite \p-energy,
which on bounded sets and for upper bounded boundary data is equivalent to
$u\in\Np(V)$.
We denote these
solutions by $\tS f$.
  Examples~6.5 and~6.6 in \cite{BBSjodin} show that $\tP f$ and $\tS f$ 
are in general different and that there are $\tP$-resolutive boundary data
which are not $\tS$-resolutive.
Corollary~7.2 in \cite{BBSjodin} shows that $\tS f$
can be equivalently defined by requiring
\eqref{eq-liminf} only for q.e.\ $x\in\bdy V$.

It follows quite easily
from the definitions that for 
upper bounded~$f$,
\begin{equation}  \label{eq-St-ge-R}
\tS f \ge \uR f.
\end{equation}
Indeed, if $u$ is admissible for $\tS f$ (and $f \not\equiv -\infty$), then 
$\ut:=\min\{u,\max_{\bdy V}  f\} \in\Np(V)$ is \p-superharmonic and
thus finely continuous in $V$
(see the proof of Proposition~\ref{prop-open-set}).
Since also
\begin{equation*} 
  \bdy_p V \subset \bdy V
  \quad \text{and} \quad 
\fineliminf_{V\ni y\to x} \ut(y) \ge \liminf_{V\ni y\to x} \ut(y),
\end{equation*}
$\ut$ is admissible for $\uR f$ according to Definition~\ref{def-P}.

In combination with
Theorem~\ref{thm-lQ<=uQ} and Proposition~\ref{prop-open-set},
inequality~\eqref{eq-St-ge-R} 
together with its analogue for lower Perron solutions implies that if a bounded
function $f$ is Sobolev-resolutive in the sense of \cite{BBSjodin},
then $f$ is resolutive for any of our four 
Perron solutions,
while our definition might give more resolutive functions.
\end{remark}


\begin{thebibliography}{99}

\bibitem{aronsson67} \art{\auth{Aronsson}{G}}
  {Extension of functions satisfying Lipschitz conditions}
  {Ark. Mat.} {6} {1967} {551--561}

\bibitem{ABcomb} \art{Bj\"orn, A.}
        {The Dirichlet problem for \p-harmonic functions on the topologist's comb}
        {Math. Z.} {279} {2015} {389--405}

\bibitem{BBbook} \book{Bj\"orn, A. \AND Bj\"orn, J.}
        {\it Nonlinear Potential Theory on Metric Spaces}
    {EMS Tracts in Mathematics {\bf 17},
        European Math. Soc., Z\"urich, 2011}

\bibitem{BBvarcap} \art{\auth{Bj\"orn}{A} \AND \auth{Bj\"orn}{J}}	
        {The variational capacity with respect to nonopen sets in metric spaces}
	{Potential Anal.} {40} {2014} {57--80}

\bibitem{BBnonopen} \art{\auth{Bj\"orn}{A} \AND \auth{Bj\"orn}{J}}	
	{Obstacle and Dirichlet problems on arbitrary nonopen sets
          in metric spaces, and fine topology}
        {Rev. Mat. Iberoam.} {31} {2015} {161--214}

\bibitem{BBLat1} \art{\auth{Bj\"orn}{A}, \auth{Bj\"orn}{J}  \AND
    \auth{Latvala}{V}}
        {The weak Cartan property for the \p-fine topology on metric spaces}
	{Indiana Univ. Math. J.} {64} {2015} {915--941}

\bibitem{BBLat3} \art{\auth{Bj\"orn}{A}, \auth{Bj\"orn}{J}  \AND
    \auth{Latvala}{V}}
        {Sobolev spaces,  fine gradients and quasicontinuity on quasiopen sets
in $\R^n$ and metric spaces}
	{Ann. Acad. Sci. Fenn. Math.} {41} {2016} {551--560}

\bibitem{BBLat2} \art{\auth{Bj\"orn}{A}, \auth{Bj\"orn}{J}  \AND
    \auth{Latvala}{V}}
        {The Cartan, Choquet and Kellogg properties of the
        fine topology on metric spaces}
	{J. Anal. Math.} {135} {2018} {59--83}

\bibitem{BBLat4} \arttoappear{\auth{Bj\"orn}{A}, \auth{Bj\"orn}{J}  \AND
    \auth{Latvala}{V}}
        {The Dirichlet problem for \p-minimizers on finely open sets in metric spaces}
        {Potential Anal.}

\bibitem{BBLat5} \artprep{\auth{Bj\"orn}{A}, \auth{Bj\"orn}{J}  \AND
    \auth{Latvala}{V}}
        {Convergence and local-to-global results for \p-superminimizers on quasiopen sets}
        {\emph{Preprint}, 2022} \\
        {\tt \href{https://arxiv.org/abs/2206.02697}{arXiv:2206.02697}}

\bibitem{BBMaly} \art{\auth{Bj\"orn}{A}, \auth{Bj\"orn}{J}  \AND
    \auth{Mal\'y}{J}}
        {Quasiopen and \p-path open sets, and characterizations of quasicontinuity}
	{Potential Anal.} {46} {2017} {181--199}

\bibitem{BBS2} \art{Bj\"orn, A., Bj\"orn, J. \AND Shan\-mu\-ga\-lin\-gam, N.}
        {The Perron method for \p-harmonic functions in metric spaces}
        {J. Differential Equations} {195} {2003} {398--429}

\bibitem{BBS5} \art{Bj\"orn, A., Bj\"orn, J. \AND Shan\-mu\-ga\-lin\-gam, N.}
        {Quasicontinuity of Newton--Sobolev functions and density of Lipschitz
        functions on metric spaces}
        {Houston J. Math.} {34} {2008} {1197--1211}

\bibitem{BBSdir} \art{Bj\"orn, A., Bj\"orn, J. \AND Shan\-mu\-ga\-lin\-gam, N.}
    {The Dirichlet problem for \p-harmonic functions with respect to
the Mazurkiewicz boundary, and new capacities}
    {J. Differential Equations} {259} {2015} {3078--3114}

\bibitem{BBSjodin} \art{\auth{Bj\"orn}{A}, \auth{Bj\"orn}{J}
      \AND \auth{Sj\"odin}{T}}
  {The Dirichlet problem for  \p-harmonic functions with respect to arbitrary
  compactifications}
  {Rev. Mat. Iberoam.} {34} {2018} {1323--1360}

 \bibitem{JB-pfine}  \art{Bj\"orn, J.} {Fine continuity on metric spaces}
         {Manuscripta Math.} {125} {2008} {369--381}

\bibitem{brelot} \art{\auth{Brelot}{M}}
        {Familles de Perron et probl\`eme de Dirichlet}
        {Acta Litt. Sci. Szeged}{9}{1939}{133--153}

\bibitem{Fug} \book{\auth{Fuglede}{B}}
         {Finely Harmonic Functions}
         {Springer, Berlin--New York, 1972}

\bibitem{GLM86} \art{Granlund, S., Lindqvist, P. \AND Martio, O.}
         {Note on the PWB-method in the nonlinear case}
         {Pacific J. Math.} {125} {1986} {381--395}

\bibitem{hansevi2}\art{\auth{Hansevi}{D}}
	{The Perron method for \p-harmonic functions
		in unbounded sets in $\R^n$ and metric spaces}
	{Math. Z.}{288}{2018}{55--74}

\bibitem{HeKiMa} \book{\auth{Heinonen}{J},
	\auth{Kilpel\"ainen}{T}
	\AND \auth{Martio}{O}}
        {Nonlinear Potential Theory of Degenerate Elliptic Equations}
        {2nd ed., Dover, Mineola, NY, 2006}

\bibitem{HKST} \book{\auth{Heinonen}{J}, \auth{Koskela}{P},
	\auth{Shan\-mu\-ga\-lin\-gam}{N} \AND \auth{Tyson}{J. T}}
       {Sobolev Spaces on Metric Measure Spaces}
	{New Mathematical Monographs {\bf 27}, Cambridge Univ. Press,
        Cambridge, 2015}

\bibitem{Kilp89} \art{Kilpel\"ainen, T.}
        {Potential theory for supersolutions of degenerate elliptic equations}
        {Indiana Univ. Math. J.} {38} {1989} {253--275}

\bibitem{KiMa92} \art{Kilpel\"ainen, T. \AND Mal\'y, J.}
        {Supersolutions to degenerate elliptic equation on quasi open sets}
        {Comm. Partial Differential Equations}
        {17} {1992} {371--405}

\bibitem{KiMa02} \art{Kinnunen, J. \AND Martio, O.}
         {Nonlinear potential theory on metric spaces}
         {Illinois Math. J.} {46} {2002} {857--883}

\bibitem{korte08} \art{Korte, R.}
         {A Caccioppoli estimate and fine continuity for superminimizers on metric spaces}
         {Ann. Acad. Sci. Fenn. Math.} {33} {2008} {597--604}

\bibitem{LatPhD} \book{\auth{Latvala}{V}}
        {Finely Superharmonic Functions of Degenerate Elliptic Equations}
        {Ann. Acad. Sci. Fenn. Ser. A I Math. Dissertationes {\bf 96}
        {(1994)}}

\bibitem{LuMaZa} \book{\auth{Luke\v{s}}{J}, \auth{Mal\'y}{J} \AND
         \auth{Zaj\'i\v{c}ek}{L}}
         {Fine Topology Methods in Real Analysis and Potential Theory}
         {Springer, Berlin--Heidelberg, 1986}

\bibitem{Perron23}\art{\auth{Perron}{O}}
	{Eine neue Behandlung der ersten Randwertaufgabe f\"ur {$\Delta u=0$}}
	{Math. Z.}{18}{1923}{42--54}

\bibitem{remak}\art{\auth{Remak}{R}}
	{\"Uber potentialkonvexe Funktionen}
	{Math. Z.}{20}{1924}{126--130}

\bibitem{Sh-harm} \art{Shan\-mu\-ga\-lin\-gam, N.}
         {Harmonic functions on metric spaces}
         {Illinois J. Math.}{45}{2001}{1021--1050}

\end{thebibliography}
\end{document}